\documentclass[11pt,english]{article}

\usepackage[numbers,comma]{natbib}
\usepackage[latin9]{inputenc}
\ifx\pdfoutput\undefined
\usepackage{graphicx}
\else
\usepackage[pdftex]{graphicx}
\fi
\usepackage{psfrag}    
\usepackage{subfigure} 
\usepackage{url}       
\usepackage{amsmath}
\usepackage{amsthm}
\usepackage{amsfonts}
\usepackage{appendix}
\usepackage{setspace}
\usepackage[top=1in, bottom=1in, left=1in, right=1in]{geometry} 
\usepackage{multirow}

\newtheorem{theorem}{Theorem}
\newtheorem{lemma}{Lemma}
\newtheorem{proposition}{Proposition}
\newtheorem{corollary}{Corollary}

\newtheorem{remark}{Remark}

\makeatletter
\usepackage{babel}
\makeatother

\begin{document}

\title{Effectively tailoring fluid and diffusion models for non-stationary state-dependent queueing systems}

\author{Young Myoung Ko and Natarajan Gautam}
\maketitle
\begin{abstract}
In this paper, we consider queueing systems where the dynamics are non-stationary and state-dependent. For performance analysis of these systems, fluid and diffusion models have been typically used. Although they are proven to be asymptotically exact, their effectiveness as approximations in the non-asymptotic regime needs to be investigated. We find that existing fluid and diffusion approximations might be either inaccurate under simplifying assumptions or computationally intractable. To address this concern, this paper focuses on developing a methodology based on adjusting the fluid model so that it provides exact mean queue lengths. Further, we provide a computationally tractable algorithm that exploits Gaussian density in order to obtain performance measures of the system. We illustrate the accuracy of our algorithm using a wide variety of numerical experiments.
\end{abstract}

\section{Introduction} \label{sec_introduction}
There are several applications of systems where the dynamics are state-dependent including the repairman problem, retrial queues, chemical reactions, epidemic models, communication networks with state-dependent routing, call centers, etc. Even assuming Markovian properties, analysis of state-dependent systems is difficult. Therefore, typically fluid and diffusion approximations are used to obtain performance measures of these systems \citep{iglehart65,ethier86,mandelbaum98-2,mandelbaum02,Whitt04,Whitt06b}. These fluid and diffusion models are obtained by utilizing Functional Strong Law of Large Numbers (FSLLN) and Functional Central Limit Theorem (FCLT) which are well summarized in \citet{bill99} and \citet{Whitt02}. Using these models, one can investigate the asymptotic behavior of the system state which can be a good approximation to the original system under certain specific conditions (e.g. heavy traffic, large number of servers, etc). The previous studies mentioned above have a common feature in that they utilize FSLLN and FCLT. They, however, have different aims, scenarios, and assumptions. This paper specifically focuses on the fluid and diffusion models for the non-stationary (i.e. time-dependent or transient) and state-dependent exponential queueing systems similar to those in \citet{mandelbaum95} and \citet{mandelbaum98}. \\
In fact, some rate functions to describe the system dynamics are of the forms, $\min(\cdot,\cdot)$ and $\max(\cdot,\cdot)$ which are not differentiable everywhere. Notice that to apply the seminal weak convergence result in \citet{kurtz78}, we require differentiability of rate functions which is not satisfied in most non-stationary and state-dependent queueing systems. To extend the theory to non-smooth rate functions, \citet{mandelbaum98} proves the weak convergence by introducing a new derivative called ``scalable Lipschtz derivative'' and provides models for several queueing systems such as Jackson networks, multiserver queues with abandonments and retrials, multiclass priority preemptive queues, etc. In addition, several sets of differential equations are also provided to obtain mean values and the covariance matrix of the limit process. It, however, turns out that the resulting sets of differential equations are not computationally tractable in general and hence the theorems cannot directly be applied to obtain performance measures numerically. In a follow-on paper, \citet{mandelbaum02} provides numerical results for queue lengths and waiting times in multiserver queues with abandonments and retrials by including an additional assumption to deal with computational tractability. Specifically, that paper assumes measure zero at a set of time points where the fluid model hits non-differentiable points, which enables applying Kurtz's diffusion models. However, as pointed out in \citet{mandelbaum02}, if the system stays in a critically loaded state (non-differentiable points) for a long time, also called ``lingering'', their approach may cause significant inaccuracy. Before describing our goal, it is worthwhile to summarize the above results and point out possible problems:
\begin{itemize}
\item On one hand, \citet{mandelbaum98} provides rigorous theory to obtain the fluid and diffusion models for the system having non-smooth rate functions. On the other hand, it is not possible to solve the resulting set of differential equations to obtain performance measures numerically.
\item The additional assumption of measure zero in \citet{mandelbaum02} (also see \citet{massey02}) provides a computationally tractable way to obtain performance measures. However, when this assumption is not valid, it might cause inaccuracy in obtaining the results.
\end{itemize}
The goal of this paper is to devise a technique that strikes a balance between accuracy and computational tractability leveraging upon the fluid and diffusion models in \citet{kurtz78}. To do so, we explain when the fluid approximations might fail to achieve accuracy from a different point of view than those considered in \citet{mandelbaum98,mandelbaum02}. Our method is irrelavent to the smoothness of rate functions and provides a condition to obtain exact estimation of mean values of the system state. We apply our methodology to several queueing systems including not only multiserver queues with abandonments and retrials considered in \citet{mandelbaum98,mandelbaum02} but also more complex queueing systems such as multiclass priority preemptive queues (considered but not numerically investigated in \citet{mandelbaum98}) and peer-based networks in multimedia distribution. Here, we emphasize that this paper does NOT aim at proving another weak convergence to a limit process but pursues providing a practically effective methodology to increase accuracy in performance measures such as mean values and the covariance matrix of the system state. Our paper has discriminating features from previous research in that we:
\begin{itemize}
	\item address possible inaccuracy of the fluid model which might occur irrespective of the smoothness of rate functions,
        \item solve the fluid model directly by providing a methodology to estimate mean values exactly unlike previous research where the fluid model is unchanged and is complemented by the expected value of the diffusion model, and
	\item devise an effective algorithm transforming the fluid and diffusion models, which achieves not only increased accuracy but also computational feasibility.
\end{itemize}
We now describe the organization of this paper. In Section \ref{sec_fluiddiffusion}, we explain the fluid and diffusion models in \citet{kurtz78} and \citet{mandelbaum98}, and describe their limited applicability in practice. In Section \ref{sec_adjustedfluid}, we provide a methodology to estimate exact mean values of the system state. However, this would not immediately result in a computationally feasible approach. For that, in Section \ref{sec_heuristic}, we explain our algorithm based on Gaussian density to achieve computational tractability and the benefits of using Gaussian density. In Section \ref{sec_application}, we show how our proposed method works for the queueing system described in \citet{mandelbaum02} by comparing our method with theirs. In Section \ref{sec_additional}, we provide some numerical results for more complex queueing systems where we have reasons to believe our methodology may not be accurate, and for these cases, we determine the performance of our approach. Finally, in Section \ref{sec_conclusion}, we make concluding remarks and explain directions for future work.
\section{Recapitulating fluid and diffusion approximations} \label{sec_fluiddiffusion}
Before explaining our results, we recapitulate fluid and diffusion approximations developed by \citet{kurtz78} that we would leverage on for our methodology. As a matter of fact, the diffusion model developed in \citet{kurtz78} is not directly applicable in many queueing systems because it requires differentiability of rate functions which is sometimes not satisfied. Therefore, we also briefly mention the result in \citet{mandelbaum98} which extends the Kurtz's result to the model involving non-smooth rate functions. Further, it is worthwhile to note that for $n\in \mathbf{N}$, the state of the queueing system $X_n(t)$ includes jumps but the limit process is continuous. Therefore, the weak convergence result that is presented is with respect to uniform topology in $D$ (\citet{bill99} and \citet{Whitt02}).\\
Let $X(t)$ be a $d$-dimensional stochastic process which is the solution to the following integral equation:
\begin{eqnarray}
	X(t) = x_0 + \sum_{i=1}^{k} l_i Y_{i} \bigg(\int_{0}^{t} f_{i}\big(s,X(s)\big)ds \bigg), \label{eqn_001}
\end{eqnarray}
where $x_0 = X(0)$ is a constant, $Y_{i}$'s are independent rate $1$ Poisson processes, $l_i \in \mathbf{Z}^d$ for $i \in \{1,2,\ldots, k\}$ are constants, and $f_i$'s are continuous functions such that $|f_i(t,x)| \le C_i (1+|x|)$ for some $C_i < \infty$, $t\le T$ and $T<\infty$.
Note that we just consider a finite number of $l_i$'s to simplify proofs, which is reasonable for real world applications.
It is usually not tractable to solve the integral equation (\ref{eqn_001}). Therefore, to approximate the $X(t)$ process, define a sequence of stochastic processes $\{X_n(t)\}$ which satisfy the following integral equation:
\begin{eqnarray*}
	X_n(t) = x_0 + \sum_{i=1}^{k} \frac{1}{n} l_i Y_{i} \bigg(\int_{0}^{t} n f_{i}\big(s, X_n(s)\big)ds \bigg). \label{eqn_002}
\end{eqnarray*}
Typically the process $X_n(t)$ (usually called a scaled process) is obtained by taking $n$ times faster rates of events and $1/n$ of the increment of the system state. This type of setting is used in the literature and is denoted as ``uniform acceleration'' in \citet{massey98}, and \citet{mandelbaum98,mandelbaum02}. Then, the following theorem provides the fluid model to which $\{X_n(t)\}$ converges almost surely as $n\rightarrow \infty$. Define
\begin{eqnarray}
	F(t,x) = \sum_{i=1}^{k} l_i f_{i}(t,x) \label{eqn_F}.
\end{eqnarray}
\begin{theorem}[Fluid model, \citet{kurtz78}] \label{theo_fluid}
If there is a constant $M < \infty$ such that $|F(t,x)-F(t,y)| \le M|x-y|$ for all $t \le T$ and $T<\infty$. Then, $\lim_{n \rightarrow \infty} X_n(t) = \bar{X}(t)$ a.s. where $\bar{X}(t)$ is the solution to the following integral equation:
\begin{eqnarray*}
	\bar{X}(t) = x_0 + \sum_{i=1}^{k} l_i \int_{0}^{t} f_{i}\big(s, \bar{X}(s)\big)ds. \label{eqn_003}
\end{eqnarray*}
\end{theorem}
Note that $\bar{X}(t)$ is a deterministic time-varying quantity. We will subsequently connect $\bar{X}(t)$ and $X(t)$ defined in equation (\ref{eqn_001}), but before that we provide the following result. Once we have the fluid model, we can obtain the diffusion model from the scaled centered process ($D_n(t)$). Define $D_n(t)$ to be $\sqrt{n}\big(X_n(t) - \bar{X}(t)\big)$. Then, the limit process of $D_n(t)$ is provided by the following theorem.
\begin{theorem}[Diffusion model, \citet{kurtz78}] \label{theo_diffusion}
If $f_i$'s and $F$, for some $M<\infty$, satisfy
\begin{eqnarray*}
	|f_i(t,x)-f_i(t,y)| \le M|x-y| \quad \textrm{and} \quad  \bigg| \frac{\partial}{\partial x_i}F(t,x)\bigg| \le M, \qquad \textrm{for } i \in \{1,\ldots,k\} \textrm{ and } 0\le t\le T,
\end{eqnarray*}
then $\lim_{n \rightarrow \infty} D_n(t) = D(t)$ where $D(t)$ is the solution to
\begin{eqnarray*}
	D(t) = \sum_{i=1}^{k} l_i \int_{0}^{t} \sqrt{f_i\big(s,\bar{X}(s)\big)}dW_i(s) + \int_{0}^{t} \partial F\big(s,\bar{X}(s)\big)D(s) ds, \label{eqn_004}
\end{eqnarray*}
$W_i(\cdot)$'s are independent standard Brownian motions, and $\partial F(t,x)$ is the gradient matrix of $F(t,x)$ with respect to $x$.
\end{theorem}
\begin{remark} \label{rem_nondiff}
Theorem \ref{theo_diffusion} requires that $F(\cdot,\cdot)$ has a continuous gradient matrix. Therefore, if we don't have such an $F$, then we cannot apply Theorem \ref{theo_diffusion} directly to obtain the diffusion model.
\end{remark}
\begin{remark} \label{rem_gaussian}
According to \citet{ethier86}, if $D(0)$ is a constant or a Gaussian random vector, then $D(t)$ is a Gaussian process.
\end{remark}
Now, we have the fluid and diffusion models for $X_n(t)$. Therefore, for a large $n$, $X_n(t)$ is approximated by
\begin{eqnarray*}
	X_n(t) \approx \bar{X}(t) + \frac{D(t)}{\sqrt{n}}. \label{eqn_005}
\end{eqnarray*}
If we follow this approximation, we can also approximate the mean and covariance matrix of $X_n(t)$ denoted by $E\big[X_n(t)\big]$ and $Cov \big[X_n(t),X_n(t)\big]$ respectively as
\begin{eqnarray}
	E\big[X_n(t)\big] &\approx& \bar{X}(t) + \frac{E\big[D(t)\big]}{\sqrt{n}}, \label{eqn_006} \\
	Cov \big[X_n(t),X_n(t)\big] &\approx& \frac{Cov \big[D(t),D(t)\big]}{n}. \label{eqn_007}
\end{eqnarray}
In equations (\ref{eqn_006}) and (\ref{eqn_007}), only $\bar{X}(t)$ is known. Therefore, in order to get approximated values of $E\big[X_n(t)\big]$ and $Cov \big[X_n(t),X_n(t)\big]$, we need to obtain $E\big[D(t)\big]$ and $Cov\big[D(t),D(t)\big]$. The following theorem provides a methodology to obtain $E\big[D(t)\big]$ and $Cov \big[D(t),D(t)\big]$.
\begin{theorem}[Mean and covariance matrix of linear stochastic systems, \citet{arnold92}] \label{theo_moment}
Let $Y(t)$ be the solution to the following linear stochastic differential equation.
\begin{eqnarray*}
	dY(t) = A(t)Y(t)dt + B(t)dW(t), \quad Y(0)=0, \label{eqn_008}
\end{eqnarray*}
where $A(t)$ is a $d \times d$ matrix, $B(t)$ is a $d \times k$ matrix, and W(t) is a $k$-dimensional standard Brownian motion.
Let $M(t) = E\big[Y(t)\big]$ and $\Sigma(t) = Cov\big[Y(t), Y(t)\big]$. Then, $M(t)$ and $\Sigma(t)$ are the solution to the following ordinary differential equations:
\begin{eqnarray}
	\frac{d}{dt}M(t) &=& A(t) M(t) \label{eqn_009} \nonumber \\
	\frac{d}{dt}\Sigma(t) &=& A(t) \Sigma(t) + \Sigma(t) A(t)' + B(t)B(t)'. \label{eqn_010}
\end{eqnarray}
\end{theorem}
\begin{corollary} \label{cor_moment}
If $M(0)=0$, then $E\big[M(t)\big] = 0$ for $t \ge 0$.
\end{corollary}
By Corollary \ref{cor_moment}, if $D(0)=0$, then $E\big[D(t)\big] = 0$ for $t \ge 0$. Therefore, if $\bar{X}(0) = X(0) = x_0$, then we can rewrite (\ref{eqn_006}) to be
\begin{eqnarray*}
	E\big[X_n(t)\big] &\approx& \bar{X}(t). \label{eqn_011}
\end{eqnarray*}
Recalling Remark \ref{rem_nondiff}, the diffusion model in \citet{kurtz78} requires differentiability of rate functions. Otherwise, we cannot apply Theorem \ref{theo_diffusion}. To address this problem, \citet{mandelbaum98} introduces a new derivative called ``scalable Lipschitz derivative'' and proves the weak convergence using it. Unlike the result in \citet{kurtz78}, it turns out that the diffusion limit may not be a Gaussian process when rate functions are not differentiable everywhere. In \citet{mandelbaum98}, expected values of the diffusion model may not be zero (compare it with Corollary \ref{cor_moment}) and could adjust the inaccuracy in the fluid model (see \citet{mandelbaum02}). The resulting differential equations for the diffusion model, however, are computationally intractable. For example, in \citet{mandelbaum98}, one of the differential equations has the following form:
\begin{eqnarray}
	\frac{d}{dt}E\big[Q_1^{(1)}(t)\big] &=& (\mu_t^1 \mathbf{1}_{\{Q_1^{(0)} \le n_t\}} + \beta_t \mathbf{1}_{\{Q_1^{(0)} > n_t\}})E\big[Q_1^{(1)}(t)^-\big] \nonumber \\ 
		&& - (\mu_t^1 \mathbf{1}_{\{Q_1^{(0)} < n_t\}} + \beta_t \mathbf{1}_{\{Q_1^{(0)} \ge n_t\}})E\big[Q_1^{(1)}(t)^+\big] + \mu_t^2E\big[Q_2^{(1)}(t)\big], \label{eqn_actdiff}
\end{eqnarray}
rendering it to be intractable.\\
Therefore, \citet{mandelbaum02}, as we understood, resorts to the method in \citet{kurtz78} by assuming measure zero at non-smooth points to avoid computational difficulty. As described in Section \ref{sec_introduction}, in this paper, our objective is to give the fluid and diffusion models a fresh look from an alternative perspective, and suitably adjust them for non-asymptotic scenarios. This is presented in the next section.
\section{Adjusted fluid model} \label{sec_adjustedfluid}
In this section, we first explain the possibility of inaccuracy when obtaining mean values of the system state using the fluid model. Then, we provide an adjusted fluid model to estimate exact mean values. For those, we consider the actual integral equation to get the exact value of $E\big[X(t)\big]$ by the following theorem.
\begin{theorem}[Expected value of $X(t)$] \label{theo_exp}
Consider $X(t)$ defined in equation (\ref{eqn_001}). Then, for $t\le T$, $E\big[X(t)\big]$ is the solution to the following integral equation.
\begin{eqnarray}
	E\big[X(t)\big] = x_0 + \sum_{i=1}^k l_i \int_{0}^{t}E\Big[f_i\big(s,X(s)\big)\Big] ds \label{eqn_013}
\end{eqnarray}
\begin{proof}
Take expectation on both sides of equation (\ref{eqn_001}). Then, 
\begin{eqnarray}
	E\big[X(t)\big] &=& x_0 + \sum_{i=1}^k l_i E\Bigg[Y_i\bigg(\int_{0}^{t}f_i\big(s,X(s)\big)ds \bigg)\Bigg] \nonumber \\
		&=& x_0 + \sum_{i=1}^k l_i E\bigg[\int_{0}^{t}f_i\big(s, X(s)\big)ds \bigg] \textrm{ since $Y_i(\cdot)$'s are nonhomogeneous Poisson processes} \nonumber \\
		&=& x_0 + \sum_{i=1}^k l_i \int_{0}^{t}E\Big[f_i\big(s, X(s)\big)\Big] ds \textrm{  by Fubini theorem in \citet{folland99}.} \nonumber
\end{eqnarray}
Therefore, we prove the theorem. 
\end{proof}
\end{theorem}
Comparing Theorems \ref{theo_fluid} and \ref{theo_exp}, notice that we cannot conclude that $\bar{X}(t)$ in Theorem \ref{theo_fluid} and $E\big[X(t)\big]$ in Theorem \ref{theo_exp} are close enough since $E\big[f_i(t, X(t))\big] \neq f_i\big(t, E[X(t)]\big)$. In some applications, $f_i$'s might be constants or linear combinations of components of $X(t)$. In those cases, Theorem \ref{theo_exp} and the following corollary imply that the fluid model would be the exact estimation of mean values of the system state.
\begin{corollary} \label{cor_exp}
If $f_i(t,x)$'s are constants or linear combinations of the components of $x$, Then,
\begin{eqnarray}
	E[X(t)] = \bar{X}(t), \nonumber
\end{eqnarray}
where $X(t)$ is the solution to (\ref{eqn_001}) and $\bar{X}(t)$ is the deterministic fluid model from theorem \ref{theo_fluid}.
\begin{proof}
Using linearity of expectation in \citet{williams91}, we can obtain the same integral equation for both $E\big[X(t)\big]$ and $\bar{X}(t)$.
\end{proof}
\end{corollary}
However, if we have different forms of $f_i$'s where $E\big[f_i(t, X(t))\big] \neq f_i\big(t, E[X(t)]\big)$, then the fluid model would be inaccurate. As seen in Section \ref{sec_fluiddiffusion}, the fluid model does not require differentiability of rate functions in both \citet{kurtz78} and \citet{mandelbaum98}. In \citet{mandelbaum98}, the diffusion model can contribute to mean values of the system state. However, as seen in equation (\ref{eqn_actdiff}), the differential equations to obtain mean values of the diffusion limit are not computationally tractable. Even if they are numerically solvable, mean values of the diffusion limit is zero by the time the fluid limit hits a non-differentiable point for the first time. We will show in Section \ref{sec_additional} that inaccuracy begins to occur before the fluid limit hits that point. Therefore, we approach this problem in a different point of view.\\ 
The basic idea of our approach is to construct a new process ($Z(t)$) so that its fluid model is exactly same as mean values of the original process $X(t)$ as described in Theorem \ref{theo_exp}. Define a set $\mathbb{F}$ of all distribution functions that have a finite mean and covariance matrix in $\mathbf{R}^d$. This set is valid for the fluid model since conditions on $f_i$'s guarantee that $E\big[|X(t)|\big] < \infty$ and $|Cov[X(t),X(t)]| < \infty$ for all $t \le T$. Define a subset $\mathbb{F}_0$ of $\mathbb{F}$  such that any $f \in \mathbb{F}_0$ has zero mean. We call an element of $\mathbb{F}_0$ a ``base distribution'' for the remainder of this paper.
\begin{proposition} \label{prop_exp}
$E\big[f_i(t,X(t))\big]$ can be represented as a function of $E[X(t)]$ for $t \le T$.
\begin{proof}
For fixed $t_0 \le T$, suppose the distribution of $X(t_0)$ is $F$. Then, $F \in \mathbb{F}$. For $F \in \mathbb{F}$, we can always find $F_0 \in \mathbb{F}_0$ such that $F(x) = F_0(x-\mu)$ where $\mu = E[X(t_0)] = \int_{\mathbf{R}^d} x dF$. Then,
\begin{eqnarray*}
	E\big[f_i(t_0,X(t_0))\big] &=& \int_{\mathbf{R}^d} f_i(t_0,x) dF \\
		&=& \int_{\mathbf{R}^d} f_i(t_0,x+\mu) dF_0.
\end{eqnarray*}
Since the integration removes $x$, by making $t_0$ and $\mu$ variables (i.e. substitute $t_0$ and $\mu$ with $t$ and $\mu(t)$ respectively), we have
\begin{eqnarray*}
	E\big[f_i(t,X(t))\big] = g_i(t,\mu(t)), \textrm{ for some function } g_i.
\end{eqnarray*}
\end{proof}
\end{proposition}
\begin{remark}
Proposition \ref{prop_exp} does not mean that $\mu(\cdot)$ completely identifies the function $g_i(\cdot,\cdot)$. In fact, the function $g_i(\cdot,\cdot)$ might be unknown unless the base distribution is identified but we can say that such a function $g_i(\cdot,\cdot)$ exists.
\end{remark}
For $t  \le T$, let $\mu(t) = E\big[X(t)\big]$. Let $g_i\big(t,\mu(t)\big) = E\big[f_i(t,X(t))\big]$ for $i \in \{1, \ldots, k\}$. Then, we can construct a new stochastic process $Z(t)$ which is the solution to the following integral equation:
\begin{eqnarray}
	Z(t) = z_0 + \sum_{i=1}^{k} l_i Y_{i} \bigg(\int_{0}^{t} g_{i}\big(s,Z(s)\big)ds \bigg). \label{eqn_014}
\end{eqnarray}
Based on equation (\ref{eqn_014}), define a sequence of stochastic processes $\{Z_n(t)\}$ satisfying
\begin{eqnarray}
	Z_n(t) = x_0 + \sum_{i=1}^{k} \frac{1}{n} l_i Y_{i} \bigg(\int_{0}^{t} n g_{i}\big(s, Z_n(s)\big)ds \bigg). \label{eqn_adjseq}
\end{eqnarray}
Next, we would like to obtain the fluid model for $Z_n(t)$. Before doing that, we, however, need to check whether the functions $g_i$'s satisfy the conditions to apply Theorem \ref{theo_fluid}. Following lemmas provide the proofs that $g_i$'s meet the conditions.
\begin{lemma} \label{lem_condition1}
If $|f_i(t,x)| \le C_i (1+|x|)$ for $t\le T$, then $g_i(t,x)$'s satisfy
\begin{eqnarray*}
	|g_i(t,x)| &\le& D_i (1+|x|) \quad \textrm{for some } D_i < \infty. \label{eqn_015} 
\end{eqnarray*}
\begin{proof}
To prove this lemma, we need to show that $E\big[|X(t)|\big] \le K\Big(1+\big|E\big[X(t)\big]\big|\Big)$ for $K < \infty$ and $t \le T$. We first show it in the one-dimensional case and then extend it to the $d$-dimensional case.\\
Let, for fixed $t_0\le T$, $X = X(t_0)$ having mean $\mu$ and variance $\sigma^2$, and $f_i(X) = f_i\big(t_0,X(t_0)\big)$. Then, by Cauchy-Schwarz inequality,
\begin{eqnarray}
	E\big[|X|\big] \le \sqrt{E[X^2]} = \sqrt{\mu^2 + \sigma^2} \le  |\mu| + \sigma \le D (1 + |\mu|) \quad \textrm{for } D = \max(1,\sigma) \label{eqn_016}.
\end{eqnarray}
Now, we have the one-dimensional case and can move to the $d$-dimensional case. Suppose $X$ has a mean vector $\mu$ and a covariance matrix $\Sigma$ such that $X=(x_1, \ldots, x_d)'$, $\mu = (\mu_1, \ldots, \mu_d)'$. Then,
\begin{eqnarray}
	E\big[|X|\big] &=& E\bigg[\sqrt{\sum_{i=1}^{d}x_i^2}\bigg] \le E\bigg[\sum_{i=1}^{d}|x_i|\bigg] = \sum_{i=1}^{d} E\big[|x_i|\big] \nonumber \\
		&\le& D\bigg(d+ \sum_{i=1}^{d} |\mu_i|\bigg) \quad \textrm{by (\ref{eqn_016})} \qquad \textrm{for } D = \max(1,\sigma_1, \ldots, \sigma_d) \nonumber \\
		&\le& D\bigg(d+ d\sqrt{\sum_{i=1}^{d} \mu_{i}^{2}}\bigg) \textrm{by Cauchy-Schwarz inequality} \nonumber \\
		&=& Dd\big(1+|\mu|\big) \label{eqn_017}.
\end{eqnarray}
Now we have $E\big[|X|\big] \le K\Big(1+\big|E[X]\big|\Big)$ for the $d$-dimensional random vector $X$ where $K=Dd$. Then,
\begin{eqnarray}
	\Big|E\big[f_i(X)\big]\Big| &\le& E\Big[\big|f_i(X)\big|\Big] \le C_i + C_i E\big[|X|\big] \quad \textrm{from assumption} \nonumber \\
		&\le& C_i + C_i K \big(1 + |\mu|\big)  \nonumber \le D_i \big(1 + |\mu|\big) \quad \textrm{for } D_i = C_i + C_iK \quad \textrm{by equation (\ref{eqn_017})} \label{eqn_018}\nonumber
\end{eqnarray}
Note $g_i(t_0,\mu) =  E\big[f_i(X)\big]$. Since $|\Sigma|$ is bounded on $t\le T$, if we make $t_0>0$ arbitrary, we prove the lemma. 
\end{proof}
\end{lemma}
For the next lemma, we would like to define
\begin{eqnarray}
  G(t,x) = \sum_{i=1}^k l_i g_i(t,x). \label{eqn_G}
\end{eqnarray}
\begin{lemma} \label{lem_condition2}
For $t\le T$, if $|f_i(t,x)-f_i(t,y)| \le M |x-y|$, then $g_i(t,x)$'s satisfy
\begin{eqnarray*}
	|g_i(t,x) - g_i(t,y)| \le M|x-y|,
\end{eqnarray*}
and if $|F(t,x)-F(t,y)| \le M |x-y|$, then $G(t,x)$ satisfies
\begin{eqnarray*}
	|G(t,x) - G(t,y)| \le M|x-y|. \label{eqn_019} 
\end{eqnarray*}
\begin{proof}
For fixed $t_0 \le T$, let $X = X(t_0)$ and $Y = Y(t_0)$ and suppose $X$ and $Y$ have a same base distribution $H_0$ (we use $H$ instead of $F$ to avoid confusion with $F$ in (\ref{eqn_F})) where $E[X] = \mu_1$ and $E[Y]=\mu_2$. Then, the distribution $H_1$ of $X$ and $H_2$ of $Y$ satisfy
\begin{eqnarray*}
	H_1(x) &=& H_0(x-\mu_1), \quad \textrm{and} \\
	H_2(y) &=& H_0(y-\mu_2),
\end{eqnarray*}
respectively.
Now, we have
\begin{eqnarray*}
	\Big|E\big[F(X)\big]-E\big[F(Y)\big]\Big| &=& \bigg|\int_{\mathbf{R}^d} F(x) dH_1 - \int_{\mathbf{R}^d} F(y) dH_2\bigg|. \label{eqn_021}
\end{eqnarray*}
By transforming variables,
\begin{eqnarray}
	\Big|E\big[F(X)\big]-E\big[F(Y)\big]\Big| &=& \bigg|\int_{\mathbf{R}^d} F(x+\mu_1) dH_0 - \int_{\mathbf{R}^d} F(y+\mu_2) dH_0\bigg| \nonumber \\
		&=& \bigg|\int_{\mathbf{R}^d} \big(F(x+\mu_1) - F(x+\mu_2)\big) dH_0\bigg| \quad \textrm{by linearity}, \nonumber \\
		&\le& \int_{\mathbf{R}^d} \bigg|\big(F(x+\mu_1) - F(x+\mu_2)\big)\bigg| dH_0 \nonumber \\
		&\le& M \int_{\mathbf{R}^d} |\mu_1-\mu_2| dH_0 = M|\mu_1 - \mu_2| \quad \textrm{by assumption}. \label{eqn_022} \nonumber
\end{eqnarray}
Note $G\big(t_0, \mu_1 \big) = E\big[F(X)\big]$ and $G\big(t_0, \mu_2 \big) = E\big[F(Y)\big]$. Then, by making $t_0>0$ arbitrary, we prove the second part, i.e. if $|F(t,x)-F(t,y)|\le M|x-y|$ then $|G(t,x)-G(t,y)|\le M|x-y|$.
We can prove the first part, i.e. if $|f_i(t,x)-f_i(t,y)| \le M |x-y|$, then $|g_i(t,x) - g_i(t,y)| \le M|x-y|$, in a similar fashion and hence we have the lemma.
\end{proof}
\end{lemma}
Lemmas \ref{lem_condition1} and \ref{lem_condition2} show that if $f_i$'s satisfy the conditions to obtain the fluid limit of $X_n(t)$, then $g_i$'s are also eligible for the fluid model of $Z_n(t)$. Therefore, we are now able to provide the adjusted fluid model based on Lemmas \ref{lem_condition1} and \ref{lem_condition2}.
\begin{theorem}[Adjusted fluid model] \label{theo_modfluid}
Assume
\begin{eqnarray}
	\big|f_i(t,x)\big| &\le& C_i\big(1+|x|\big) \quad \textrm{for } i\in \{1,\ldots, k\}, \label{eqn_024}\\
	\big|F(t,x)-F(t,y)\big| &\le& M|x-y|. \label{eqn_025}
\end{eqnarray}
Then, $\lim_{n \rightarrow \infty} Z_n(t) = \bar{Z}(t)$ a.s., where $\bar{Z}(t)$ is the solution to the following integral equation:
\begin{eqnarray}
	\bar{Z}(t) = x_0 + \sum_{i=1}^{k} l_i \int_{0}^{t} g_{i}\big(s, \bar{Z}(s)\big)ds, \label{eqn_026}
\end{eqnarray}
and furthermore
\begin{eqnarray}
	\bar{Z}(t) = E\big[X(t)\big] = x_0 + \sum_{i=1}^k l_i \int_{0}^{t}E\Big[f_i\big(s,X(s)\big)\Big] ds. \label{eqn_027}
\end{eqnarray}
\begin{proof}
From Lemmas \ref{lem_condition1} and \ref{lem_condition2}, (\ref{eqn_024}) and (\ref{eqn_025}) imply
\begin{eqnarray*}
	|g_i(t,x)|\le D_i (1+|x|) \quad \textrm{and} \quad |G(t,x) - G(t,y)| \le M|x-y|. \label{eqn_029} 
\end{eqnarray*}
Therefore, by Theorem \ref{theo_fluid}, we have equation (\ref{eqn_026}), and by definition of $g_i(t,x)$'s, we have equation (\ref{eqn_027}).
\end{proof}
\end{theorem}
In Theorem \ref{theo_modfluid}, we have the same conditions for $f_i$'s and $g_i$'s, and $g_i(t,x)$'s do not guarantee $E\big[g_i(t,X(t))\big] = g_i\big(t,E[X(t)]\big)$ either. However, comparing equation (\ref{eqn_027}) with equation (\ref{eqn_013}) in Theorem \ref{theo_exp}, we notice that Theorem \ref{theo_modfluid} via equation (\ref{eqn_027}) could provide the exact estimation of $E\big[X(t)\big]$.\\
Though Theorem \ref{theo_modfluid} provides the exact estimation of $E\big[X(t)\big]$, we should identify the functions $g_i$'s in order to obtain these values numerically or analytically. The $g_i$'s, however, cannot be identified unless the base distribution is known, which forces us to develop an algorithm to find $g_i$'s. The following section will describe our Gaussian-based method which would also be useful to adjust the diffusion model.  
\section{Gaussian-based adjustment} \label{sec_heuristic}
In Theorem \ref{theo_modfluid}, we encounter a fundamental problem in finding $g_i$'s, i.e. we need to characterize the distribution of $X(t)$. There, however, is no clear way to find the exact distribution of $X(t)$ in general. Therefore, our proposed method starts with assuming the distribution of $X(t)$. Recall that in Section \ref{sec_fluiddiffusion}, $X_n(t)$ is approximated by a Gaussian process when $x_0$ is a constant for a large $n$. Though it is not true for $n=1$, we use a Gaussian density function to obtain $g_i$'s with $z_0 = x_0$ since using the Gaussian density function provides following three benefits:
\begin{enumerate}
	\item In many applications, empirical densities are close to Gaussian density even if rate functions are not differentiable (see \citet{mandelbaum98-2}, \citet{mandelbaum02}).
	\item Gaussian distribution can be completely characterized by the mean and covariance matrix which can be obtained from the fluid and diffusion models.
	\item By using Gaussian density, $g_i$'s can achieve smoothness even if $f_i$'s are not smooth, which enable us to apply Theorem \ref{theo_diffusion} directly.
\end{enumerate}
The third benefit is not obvious and hence we provide the proof of that.
\begin{lemma} \label{lem_smooth}
Let $g_i$'s be the rate functions of $Z(t)$ obtained from Gaussian density. Then, $g_i$'s are differentiable everywhere.
\begin{proof}
Define
\begin{eqnarray}
	\phi(x,y) = \frac{1}{(2\pi)^{n/2}|\Sigma|^{1/2}}\exp \bigg(-\frac{(y-x)'\Sigma^{-1}(y-x)}{2}\bigg).\nonumber 
\end{eqnarray}
Using Gaussian density,
\begin{eqnarray*}
	g_i(t,x) = \int_{\mathbf{R}^d} f_i(t,y) \phi(x,y) dy.
\end{eqnarray*}
For $j \in \{1,\ldots, d\}$, since $\phi(x,y)$ is differentiable with respect to $x_j$ and $|f_i(t,y) \frac{d}{dx_j} \phi(x,y)|$ is integrable,
\begin{eqnarray}
	\frac{d}{dx_j}g_i(t,x) &=& \frac{d}{dx_j}\int_{\mathbf{R}^d} f_i(t,y) \phi(x,y) dy \nonumber\\
		&=& \int_{\mathbf{R}^d} f_i(t,y) \frac{d}{dx_j} \phi(x,y) dy \quad \textrm{by applying Theorem 2.27 in \citet{folland99}}, \label{eqn_folland}
\end{eqnarray}
where $x_j$ is $j^{\textrm{th}}$ component of $x$.\\
Therefore, $g_i$ is differentiable with respect to $x_j$. \\
\end{proof}
\end{lemma}
Now, we have $g_i(\cdot,\cdot)$'s which are differentiable. Then, we can apply Theorem \ref{theo_diffusion} to obtain the diffusion model for $Z_n(t)$.
\begin{proposition}[Adjusted diffusion model] \label{prop_moddiffusion} Let $g_i(\cdot,\cdot)$'s be the rate functions in $Z(t)$ obtained from Gaussian density. Define a sequence of scaled centered processes $\{V_n(t)\}$ for $t \le T$ to be
\begin{eqnarray*}
	V_n(t) = \sqrt{n}\big(Z_n(t)-\bar{Z}(t)\big), \label{eqn_030}
\end{eqnarray*}
where $Z_n(t)$ and $\bar{Z}(t)$ are solutions to equations (\ref{eqn_adjseq}) and (\ref{eqn_026}) respectively. If $f_i(t,x)$'s and $F(t, x)$ satisfy equations (\ref{eqn_024}) and (\ref{eqn_025}) respectively, then
$\lim_{n \rightarrow \infty} V_n(t) = V(t)$, where
\begin{eqnarray*}
	V(t) = \sum_{i=1}^{k} l_i \int_{0}^{t} \sqrt{g_i\big(s,\bar{Z}(s)\big)}dW_i(s) + \int_{0}^{t} \partial G\big(s,\bar{Z}(s)\big) ds, \label{eqn_031}
\end{eqnarray*}
$W_i(\cdot)$'s are independent standard Brownian motions, and $\partial G\big(t,\bar{Z}(t)\big)$ is the gradient matrix of $G\big(t,\bar{Z}(t)\big)$ with respect to $\bar{Z}(t)$. Furthermore, $V(t)$ is a Gaussian process.
\begin{proof}
From definition of $G(t,x)$ in (\ref{eqn_G}), we can easily verify that $G(t,x)$ is differentiable by Lemma \ref{lem_smooth} and hence $|G(t,x) - G(t,y)| \le M|x-y|$ implies
\begin{eqnarray*}
	\bigg|\frac{\partial}{\partial x_i} G(t,x) \bigg| \le M_i \quad \textrm{for some } M_i < \infty, t \le T, \textrm{ and } i\in \{1,\ldots, d\}.
\end{eqnarray*}
Therefore, by Theorem \ref{theo_diffusion}, we prove this proposition.
\end{proof}
\end{proposition}
\begin{corollary} \label{cor_samedistribution}
If $f_i$'s are constants or linear combinations of the components of $X(t)$. Then,
\begin{eqnarray*}
	X(t) = Z(t) \quad \textrm{in distribution}. \label{eqn_032}
\end{eqnarray*}
\begin{proof}
Using the linearity of expectation, we can verify $g_i(t,x)= f_i(t,x)$ for $i\in\{1,\ldots,k\}$.
\end{proof}
\end{corollary}
Now we have the adjusted fluid and diffusion models by utilizing Gaussian density. Therefore, instead of assuming measure zero at a set of non-differentiable points (as done in \citet{mandelbaum02}), we compare the adjusted models with the empirical mean and covariance matrix. Note when we explain Theorem \ref{theo_modfluid}, we do not consider $\Sigma(t)$, the covariance matrix of $X(t)$. However, from Gaussian density, we know that $\Sigma(t)$ characterizes the base distribution and it can be obtained from Proposition \ref{prop_moddiffusion}. Therefore, we rewrite $g_i$'s to be functions of $t$, $\bar{Z}(t)$, and $\Sigma(t)$; i.e.
\begin{eqnarray}
	g_{i}\big(t, \bar{Z}(t)\big) &\rightarrow& g_{i}\big(t, \bar{Z}(t), \Sigma(t) \big) \quad \textrm{for } i\in\{1,\ldots,k\} \textrm{ and}  \label{eqn_036}\\
	G\big(t, \bar{Z}(t)\big) &\rightarrow& G\big(t, \bar{Z}(t), \Sigma(t) \big). \label{eqn_037}
\end{eqnarray}
\begin{proposition}[Mean and covariance matrix] \label{prop_modmoment}
Let $Y(t) = \bar{Z}(t) + V(t)$. Then,
\begin{eqnarray}
	E\big(Y(t)\big) &=& \bar{Z}(t) \quad \textrm{and} \label{eqn_038} \\
	Cov\big(Y(t),Y(t)\big) &=& Cov\big(V(t),V(t)\big) = \Sigma(t). \label{eqn_039}
\end{eqnarray}
The quantities $\bar{Z}(t)$ and $\Sigma(t)$ are obtained by solving the following simultaneous ordinary differential equations with initial values given by $\bar{Z}(0) = x_0$ and $\Sigma(0)=0$:
\begin{eqnarray}
	\frac{d}{dt}\bar{Z}(t) &=& \sum_{i=1}^{k} l_i g_{i}\big(t, \bar{Z}(t), \Sigma(t) \big), \label{eqn_040} \\
	\frac{d}{dt}\Sigma(t) &=& A(t) \Sigma(t) + \Sigma(t) A(t)' + B(t)B(t)', \label{eqn_041}
\end{eqnarray}
where $A(t)$ is the gradient matrix of $G\big(t,\bar{Z}(t), \Sigma(t)\big)$ with respect to $\bar{Z}(t)$, and $B(t)$ is the $d \times k$ matrix such that its $i^\textrm{th}$ column is $l_i \sqrt{g_i\big(t,\bar{Z}(t), \Sigma(t)\big)}$.
\begin{proof}
Since $V(0) = 0$, from Corollary \ref{cor_moment}, we have (\ref{eqn_038}) and (\ref{eqn_039}). By rewriting (\ref{eqn_026}) in Theorem \ref{theo_modfluid} as a differential equation form, we have (\ref{eqn_040}), and by Theorem \ref{theo_moment}, we have (\ref{eqn_041}). Note that since both $\bar{Z}(t)$ and $\Sigma(t)$ are variables, we should solve (\ref{eqn_040}) and (\ref{eqn_041}) simultaneously.
\end{proof}
\end{proposition}
In conclusion, we define an adjusted process $Z(t)$ in Section \ref{sec_adjustedfluid} to obtain the exact $E\big[X(t)\big]$ for $X(t)$ process which is the state of a non-stationary and state-dependent queueing system. It, however, is not possible to obtain such $g_i$'s and hence in this section, we provided an algorithm by utilizing Gaussian density. From this, the limit process turns out to be a Gaussian process. We recognize that this is not true for the original process. As mentioned in Section \ref{sec_introduction}, however, this paper does not pursue finding the exact distribution of the original process but proposes an effective way to estimate mean values and the covariance matrix of the original process. Therefore, in the following sections, by means of numerical examples, we illustrate our methodology and show its effectiveness.
\section{Multiserver queues with abandonments and retrials} \label{sec_application}
\begin{figure}[htbp]
  \begin{center}
      \subfigure[ Mean numbers by assuming measure zero]{
		\includegraphics[width = .48\textwidth]{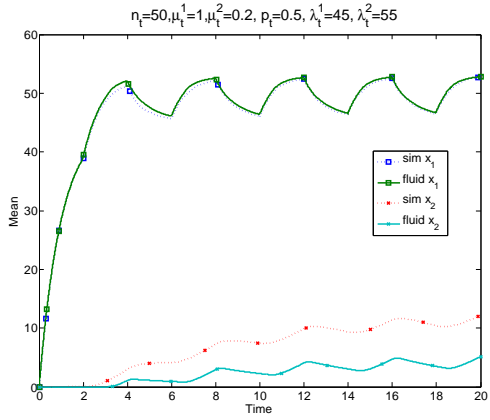}}
      \subfigure[ Mean numbers by our proposed method]{
		\includegraphics[width = .48\textwidth]{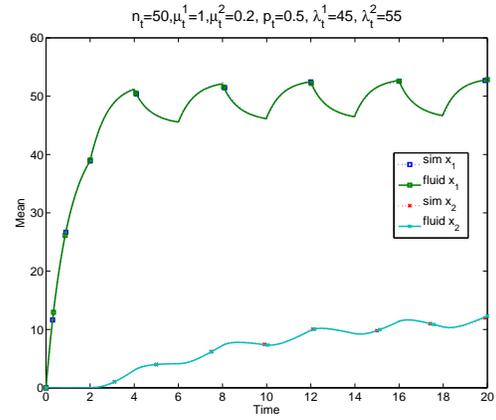}}
    \caption{Comparison of mean values, $E\big[X(t)\big]$}
    \label{fig_retrial_mean}
  \end{center}
\end{figure}

\begin{figure}[htbp]
  \begin{center}
      \subfigure[ Covariance matrix by assuming measure zero]{
		\includegraphics[width = .48\textwidth]{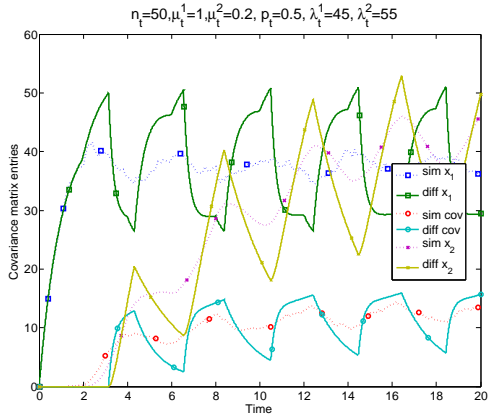}}
      \subfigure[ Covariance matrix by our proposed method]{
		\includegraphics[width = .48\textwidth]{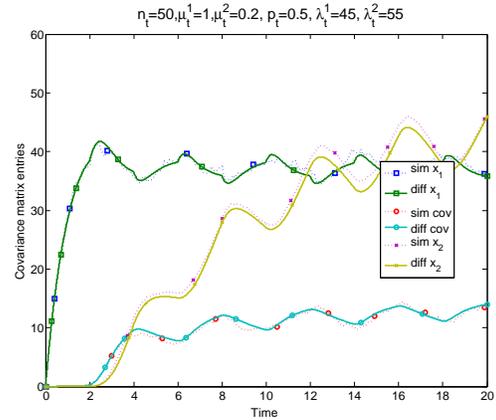}}
    \caption{Comparison of covariance matrix entries, $Cov\big[X(t),X(t)\big]$}
    \label{fig_retrial_covariance}
  \end{center}
\end{figure}

Multiserver queues with abandonments and retrials are extensively studied in the literature since they are used to model an important application, namely ``call centers'' (e.g. \citep{Halfin81,Garnet02,Whitt04,zeltyn05,Whitt06b}). In this section, therefore, we provide in-depth explanation of how our approach works in this queueing system by numerical examples.
\begin{figure}
\centering
\includegraphics[width = .6\textwidth]{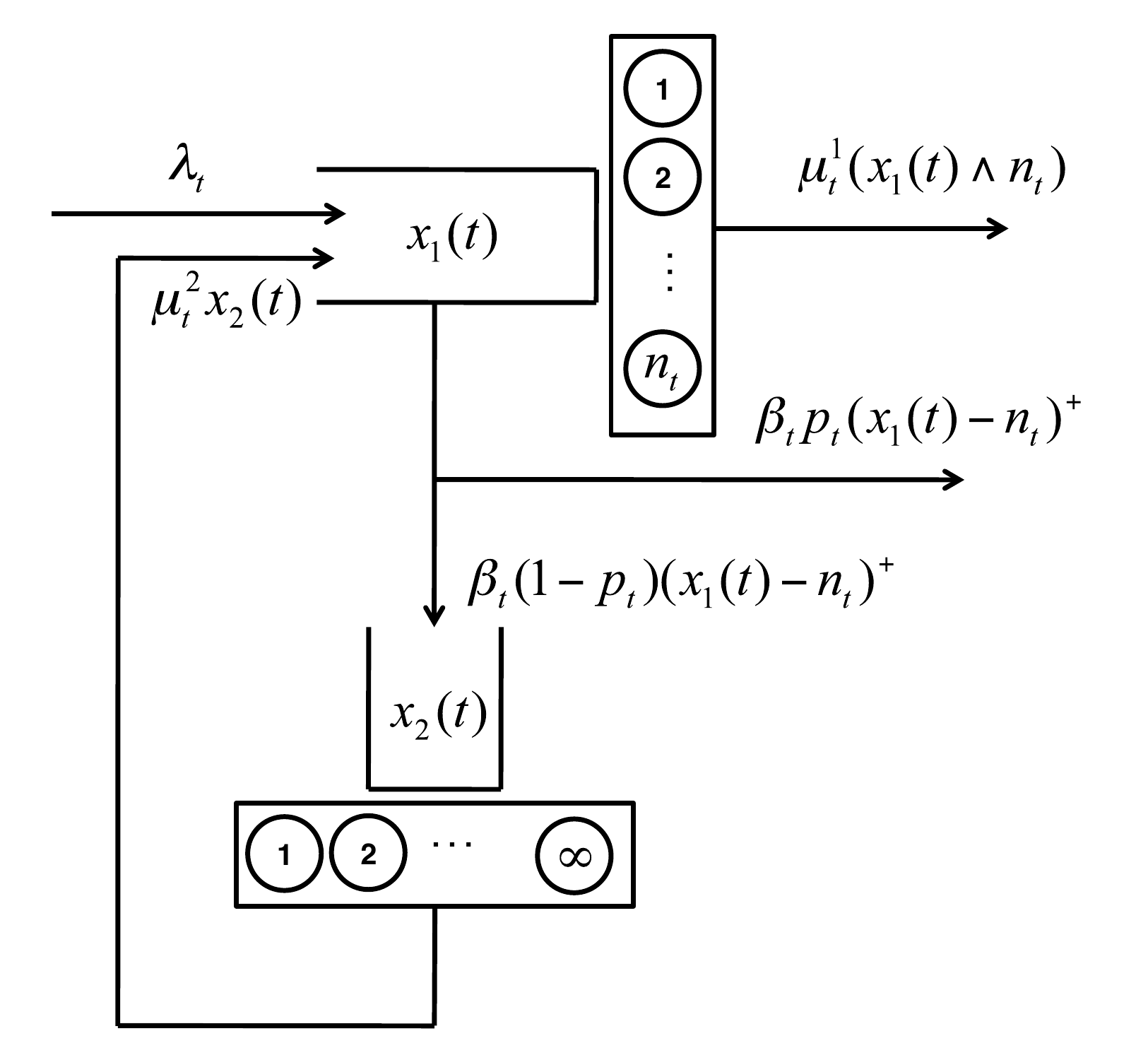}
\caption{Multiserver queue with abandonment and retrials, \citet{mandelbaum02}} \label{fig_retrial}
\end{figure}
Figure \ref{fig_retrial} illustrates a multiserver queue with abandonments and retrials described in \citet{mandelbaum98, mandelbaum02}. There are $n_t$ number of servers in the service node at time $t$. Customers arrive to the service node according to a nonhomogeneous Poisson process at rate $\lambda_t$. The service time of each customer follows a distribution having a memoryless property at rate $\mu_t^1$. Customers in the queue are served under the FCFS policy and the abandonment rate of customers is $\beta_t$ with exponentially distributed time to abandon. Abandoning customers leave the system with probability $p_t$ or go to the retrial queue with probability $1-p_t$. The retrial queue is an infinite-server-queue and hence each customer in the retrial queue waits there for a random amount of time with mean $1/\mu_t^2$ and returns to the service node.\\
Let $X(t) = \big(x_1(t),x_2(t)\big)$ be the system state where $x_1(t)$ is the number of customers in the service node and $x_2(t)$ is the number of customers in the retrial queue. Then, $X(t)$ is the unique solution to the following integral equation:
\begin{eqnarray*}
  x_1(t) &=& x_1(0) +Y_1\Big(\int_{0}^{t}\lambda_s ds\Big) + Y_2\Big(\int_{0}^{t}x_2(s)\mu_s^2ds\Big) - Y_3\Big(\int_{0}^{t}\big(x_1(s)\wedge n_s\big)\mu_s^1ds\Big) \nonumber \\
  && - Y_4\Big(\int_{0}^{t}\big(x_1(s)-n_s\big)^+\beta_s(1-p_s)ds\Big) - Y_5\Big(\int_{0}^{t}\big(x_1(s)-n_s\big)^+\beta_sp_sds\Big), \label{eqn_rx1} \\
  x_2(t) &=& x_2(0) + Y_4\Big(\int_{0}^{t}\big(x_1(s)-n_s\big)^+\beta_s(1-p_s)ds\Big) - Y_2\Big(\int_{0}^{t}x_2(s)\mu_s^2ds\Big). \label{eqn_rx2}
\end{eqnarray*}
Then, following the notation in Section \ref{sec_fluiddiffusion}, we have, for $x=(x_1,x_2)$ and $t\le T$, 
\begin{eqnarray}
  f_1(t,x) &=& \lambda_t, \nonumber \\
  f_2(t,x) &=& \mu_t^2 x_2, \nonumber \\
  f_3(t,x) &=& \mu_t^1(x_1 \wedge n_t), \nonumber \\
  f_4(t,x) &=& \beta_t(1-p_t)(x_1-n_t)^+, \quad \textrm{and} \nonumber \\
  f_5(t,x) &=& \beta_tp_t(x_1-n_t)^+. \nonumber
\end{eqnarray}
We can verify that all $f_i$'s satisfy the conditions to apply Theorem \ref{theo_fluid}. However, we cannot apply Theorem \ref{theo_diffusion} directly since $f_3$, $f_4$, and $f_5$ are not differentiable at $x_1=n_t$. To resolve this, \citet{mandelbaum02} assumes measure zero at a set of time points when the fluid limit hits the non-differentiable points and apply Theorem \ref{theo_diffusion}. \citet{mandelbaum02} addresses that for a system of a fixed size, assuming measure zero works well when $x_1(t)$ does not stay too long near the critically loaded phase. It also provides the actual form of differential equations for the diffusion model in \citet{mandelbaum98} which, in fact, are not computationally tractable, e.g. see equations (4.1) and (4.2) in \citep{mandelbaum02}. As mentioned in Section \ref{sec_adjustedfluid}, we approch the problem under a different point of view. Notice that in addition to their non-differentiability, $f_3$, $f_4$, and $f_5$ do not satisfy $E\big[f_i\big(t,X(t)\big)\big] = f_i\big(t,E\big[X(t)\big]\big)$ either. Therefore, we would like to apply Theorem \ref{theo_modfluid} to obtain $E\big[X(t)\big]$ exactly. Recalling Section \ref{sec_heuristic}, however, obtaining exact $g_i$'s is not possible and hence we obtain $g_i$'s from Gaussian density as follows:
\begin{eqnarray}
  g_1(t,x) &=& \lambda_t, \nonumber \\
  g_2(t,x) &=& \mu_t^2 x_2, \nonumber \\
  g_3(t,x) &=& \mu_t^1\big(n_t + (x_1-n_t)\Phi(n_t,x_1,\sigma_{1_t}) - \sigma_{1_t}^2 \phi(n_t,x_1,\sigma_{1_t})\big), \nonumber \\
  g_4(t,x) &=& \beta_t(1-p_t)\Big((x_1-n_t)\big(1-\Phi(n_t,x_1,\sigma_{1_t})\big)+\sigma_{1_t}^2\phi(n_t,x_1,\sigma_{1_t})\Big), \quad \textrm{and} \nonumber \\
  g_5(t,x) &=& \beta_tp_t\Big((x_1-n_t)\big(1-\Phi(n_t,x_1,\sigma_{1_t})\big)+\sigma_{1_t}^2\phi(n_t,x_1,\sigma_{1_t})\Big), \nonumber
\end{eqnarray}
where $\Phi(a,b,c)$ and $\phi(a,b,c)$ are function values at point $a$ of the Gaussian CDF and PDF respectively with mean $b$ and standard deviation $c$.\\
Since $f_1(t,x)$ and $f_2(t,x)$ are constant and linear with respect to $x$ respectively,  $g_1(t,x)=f_1(t,x)$ and $g_2(t,x)=f_2(t,x)$. The derivation of other $g_i(\cdot,\cdot)$'s is straightforward but requires some computational efforts and hence we provide the details in Appendix \ref{app_gi}. Note $g_3$, $g_4$, and $g_5$ include $\sigma_{1_t}$ which is currently treated as a function of $t$ but is used by the adjusted diffusion model (see equations (\ref{eqn_036}) and (\ref{eqn_037})). \\
For a fixed size of the system, both our proposed method and the method assuming measure zero do not guarantee the exact estimation of the system state. However, these two methods provide computational tractability. Therefore, we compare our method against the method assuming zero in \citet{mandelbaum02}. We conduct simulations under the similar settings in \citet{mandelbaum02}. We use 5,000 independent simulation runs and compare the simulation result with both our method and the method assuming measure zero. We use the constant rates for the parameters except the arrival rate. The arrival rate alternates between $45$ and $55$ every two time units. Figures \ref{fig_retrial_mean} and \ref{fig_retrial_covariance} show the estimation of mean values from one experiment. The number of servers ($n_t$) is $50$ and the service rate of each server is $1$.
As seen in Figure \ref{fig_retrial_mean}, the number of customers in service node ($x_1(t)$) stays near the critically loaded point for a long time. As \citet{mandelbaum02} points out, the method assuming measure zero shows significant difference in estimating $E\big[x_2(t)\big]$. On the other hand, our proposed method provides accurate results. When we see the estimation of the covariance matrix, we notice the similar results as the estimation of mean values. 
As seen in Figure \ref{fig_retrial_covariance}, the method assuming measure zero causes ``spikes'' as also pointed out in \citet{mandelbaum02}. Our proposed method, however, provides reasonable accuracy and no spikes at all. To verify the effectiveness of our method, we conduct several experiments with different parameter combinations.
\begin{table}[htdp]
\caption{Experiments setting}
\begin{center}
\begin{tabular}{|c|c|c|c|c|c|c|c|c|c|c|}
\hline
exp \# & svrs & $\lambda_1$ & $\lambda_2$ & $\mu_1$ & $\mu_2$ & $\beta$ & $p$ & alter & time \\ \hline \hline
1 & 50 & 40 & 80 & 1 & 0.2 & 2.0 & 0.5 & 2 & 20 \\ \hline
2 & 50 & 40 & 60 & 1 & 0.2 & 2.0 & 0.5 & 2 & 20 \\ \hline
3 & 100 & 80 & 120 & 1 & 0.2 & 2.0 & 0.7 & 2 & 20 \\ \hline
4 & 100 & 90 & 110 & 1 & 0.2 & 2.0 & 0.7 & 2 & 20 \\ \hline
5 & 50 & 40 & 80 & 1 & 0.2 & 1.5 & 0.7 & 2 & 20 \\ \hline
6 & 50 & 40 & 60 & 1 & 0.2 & 1.5 & 0.7 & 2 & 20 \\ \hline
7 & 50 & 45 & 55 & 1 & 0.2 & 2.0 & 0.5 & 2 & 20 \\ \hline
8 & 100 & 95 & 105 & 1 & 0.2 & 2.0 & 0.5 & 2 & 20 \\ \hline
9 & 150 & 140 & 160 & 1 & 0.2 & 2.0 & 0.5 & 2 & 20 \\ \hline
10 & 150 & 100 & 190 & 1 & 0.2 & 2.0 & 0.5 & 2 & 20 \\
\hline
\end{tabular}
\end{center}
\label{tab_settings}
\end{table}
Table \ref{tab_settings} describes the setting of each experiment. In Table \ref{tab_settings}, ``svrs'' is the number of servers ($n_t$), ``alter'' is the time length for which each arrival rate lasts, and ``time'' is the end time of our analysis. We already recognize that the method assuming measure zero works well when it \emph{does not linger} too long near the non-differentiable points. For comparison, therefore, our experiments contain several cases where the system \emph{does linger} relatively long around those points as well as the cases where it does not. Experiments 1-4 are intended to see the effects of ``lingering'' around non-differentiable points. We change $\beta_t=\beta$ and $p_t=p$ as well as the arrival rates in experiments 5-8 to see the effects of other parameters. In fact, from the other experements not listed in Table \ref{tab_settings}, it turns out that changing other parameters does not affect estimation accuracy significantly. Experiments 9 and 10 are set to observe how larger arrival rates and number of servers affect estimation accuracy along with the lingering around the non-differentiable points by increasing both of them. Here we explain the overall results: for the details of numerical results, see Table \ref{tab_mean_x1}-\ref{tab_var_x2} in Appendix \ref{app_table}. Similar to the results in Figures \ref{fig_retrial_mean} and \ref{fig_retrial_covariance}, we observe that lingering does debase the quality of approximations significantly when assuming measure zero. On the other hand, we see that our proposed method provides excellent accuracy in both mean and covariance values. Even if we increase both arrival rates and number of servers, we notice that lingering still affects estimation accuracy significantly when assuming measure zero but it does not in our proposed method.
\begin{figure}
\centering
\includegraphics[width = .8\textwidth]{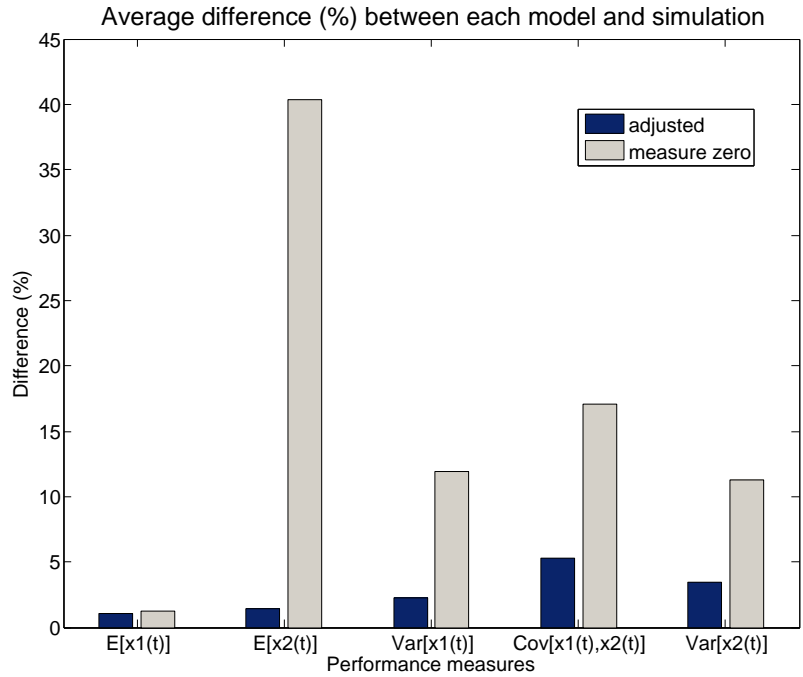}
\caption{Average difference against simulation} \label{fig_avgresult}
\end{figure}
Figure \ref{fig_avgresult} illustrates the average percentile difference of both methods against the simulation. The average difference are obtained by averaging all differences in the tables, so it does not provide the absolute comparison between two methods. However, from Figure \ref{fig_avgresult}, we intuitively notice that our proposed method shows great effectiveness relative to the method assuming measure zero.

\section{Additional Applications} \label{sec_additional}
We applied our proposed method to a wide variety of non-stationary state-dependent queueing systems. Since our proposed method is based on the adjusted fluid model, we observe that the mean queue lengths are accurate in all the systems. Also, when the rate functions are smooth or Gaussian density approximation is accurate or both, our adjusted diffusion model also provides accurate results. Due to space restrictions and our perception of how much value those cases would add, we have omitted presenting them here. Instead, we focus on scenarios with non-smooth rate functions where we conjecture Gaussian density would perhaps be inaccurate. We specifically consider two such applications \emph{not to showcase the effectiveness of our methodology, but to illustrate that there is room for improvement} for researchers in future to consider. We would like to nonetheless point out that to the best of our knowledge our approximations are still more accurate than those in the present literature. In particular, we consider multiclass preemptive priority queues (Section \ref{subsec_priority}) and peer networks (Section \ref{subsec_p2p}). For these applications, we do not provide as much detail as the multiserver queues with abandonments and retrials in Section \ref{sec_application} and show results graphically.
\subsection{Multiclass preemptive priority queues} \label{subsec_priority}
\begin{figure}
\centering
\includegraphics[width = .8\textwidth]{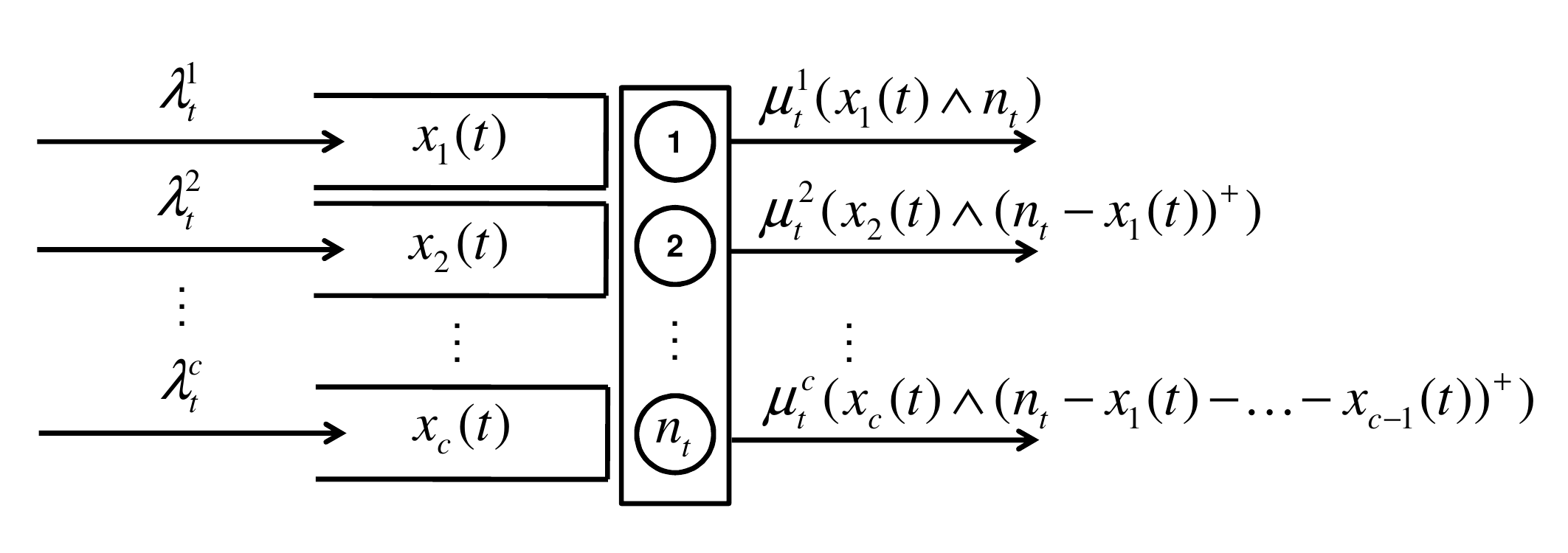}
\caption{Multiclass priority queue with preemptive policy, \citet{mandelbaum98}} \label{fig_multiclass}
\end{figure}

In this section, we consider a multiclass priority queue (see Figure \ref{fig_multiclass}) with preemptive policy which in fact is described in \citet{mandelbaum98}. It is crucial to notice that \citet{mandelbaum98} does not numerically solve this example. However, we will not only use our methodology but also extend the method assuming measure zero in \citet{mandelbaum02} for this case.\\
Explaining the priority queue briefly, there are $c$ number of classes of customers. The class $i$ customers arrive to the system with rate $\lambda_t^i$ and are served by available servers among $n_t$ number of servers with rate $\mu_t^i$ at time $t$. If a class $i$ customer arrives and there is no available server, then the highest class customer (i.e. lowest priority) in service is pushed back to the queue and the class $i$ customer is served. If there is no higher class customer, then the class $i$ customer waits in the queue. \\
In our numerical study, we use two classes of customers. Let $X(t)=\big(x_1(t),x_2(t)\big)$ be the state of the system at time $t$ where $x_1(t)$ and $x_2(t)$ are the number of customers of class 1 and 2 respectively. Then, $X(t)$ is the solution to the following integral equations:
\begin{eqnarray*}
  x_1(t) &=& x_1(0) + Y_1\Big(\int_0^t \lambda_s^1 ds\Big) - Y_3\Big(\int_0^t \mu_s^1\big(x_1(s) \wedge n_s)ds\Big), \\
  x_2(t) &=& x_2(0) + Y_2\Big(\int_0^t \lambda_s^2 ds\Big) - Y_4\Big(\int_0^t \mu_s^2\big(x_2(s) \wedge (n_s-x_1(s))^+\big)ds\Big).
\end{eqnarray*}
We set the number of servers ($n_t$) to be $200$. The arrival rate of class 1 customers ($\lambda_t^1$) is alternating between $120$ and $200$ every two time units and the arrival rate of class 2 customers ($\lambda_t^2)$ is $20$. The service rates of both class customers are $1$, i.e. $\mu_t^1 = \mu_t^2 = 1$. We conduct 5,000 simulation runs and obtain mean values by averaging them.
\begin{figure}[htbp]
  \begin{center}
      \subfigure[ Mean numbers by assuming measure zero]{
		\includegraphics[width = .48\textwidth]{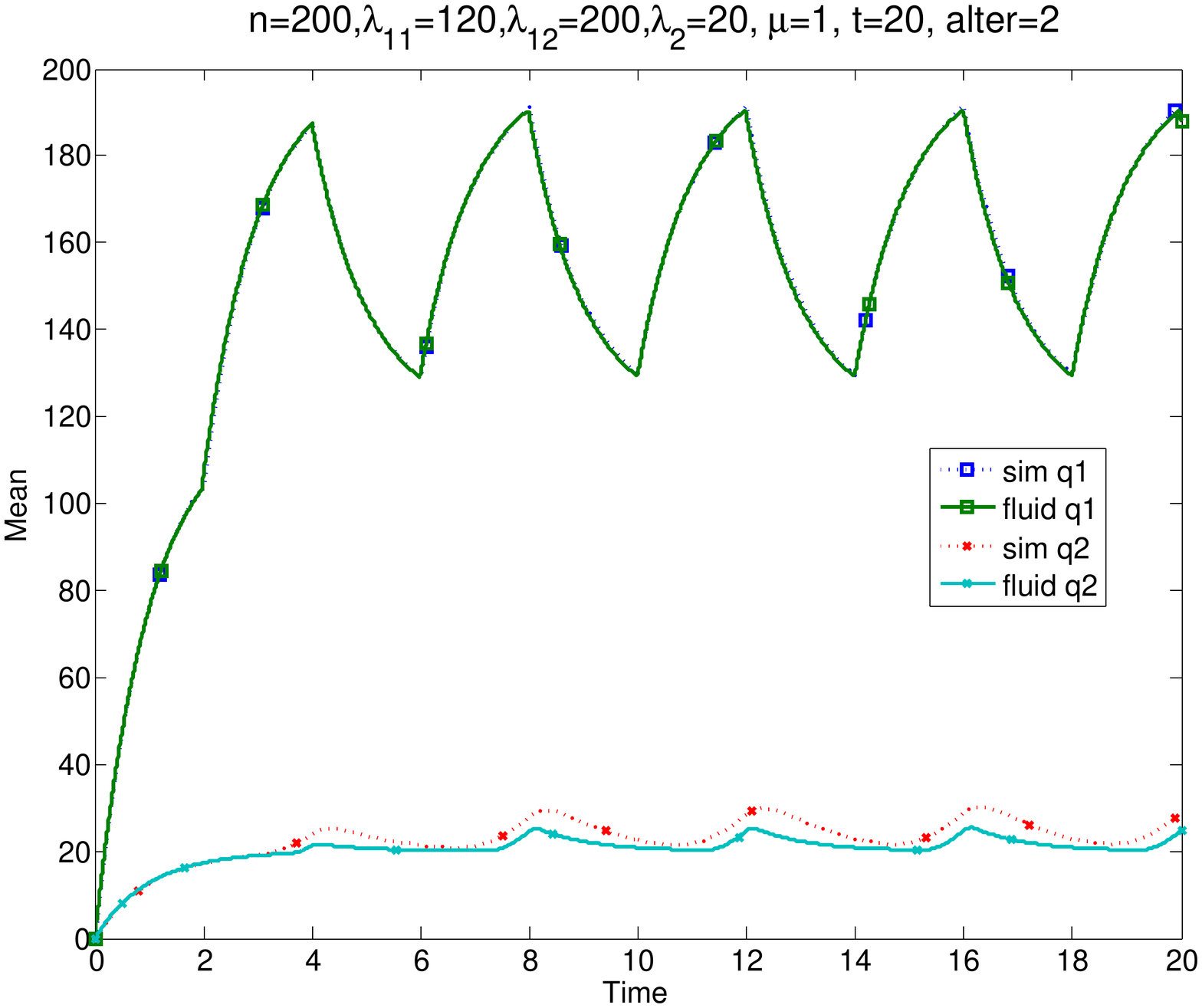}}
      \subfigure[ Mean numbers by our proposed method]{
		\includegraphics[width = .48\textwidth]{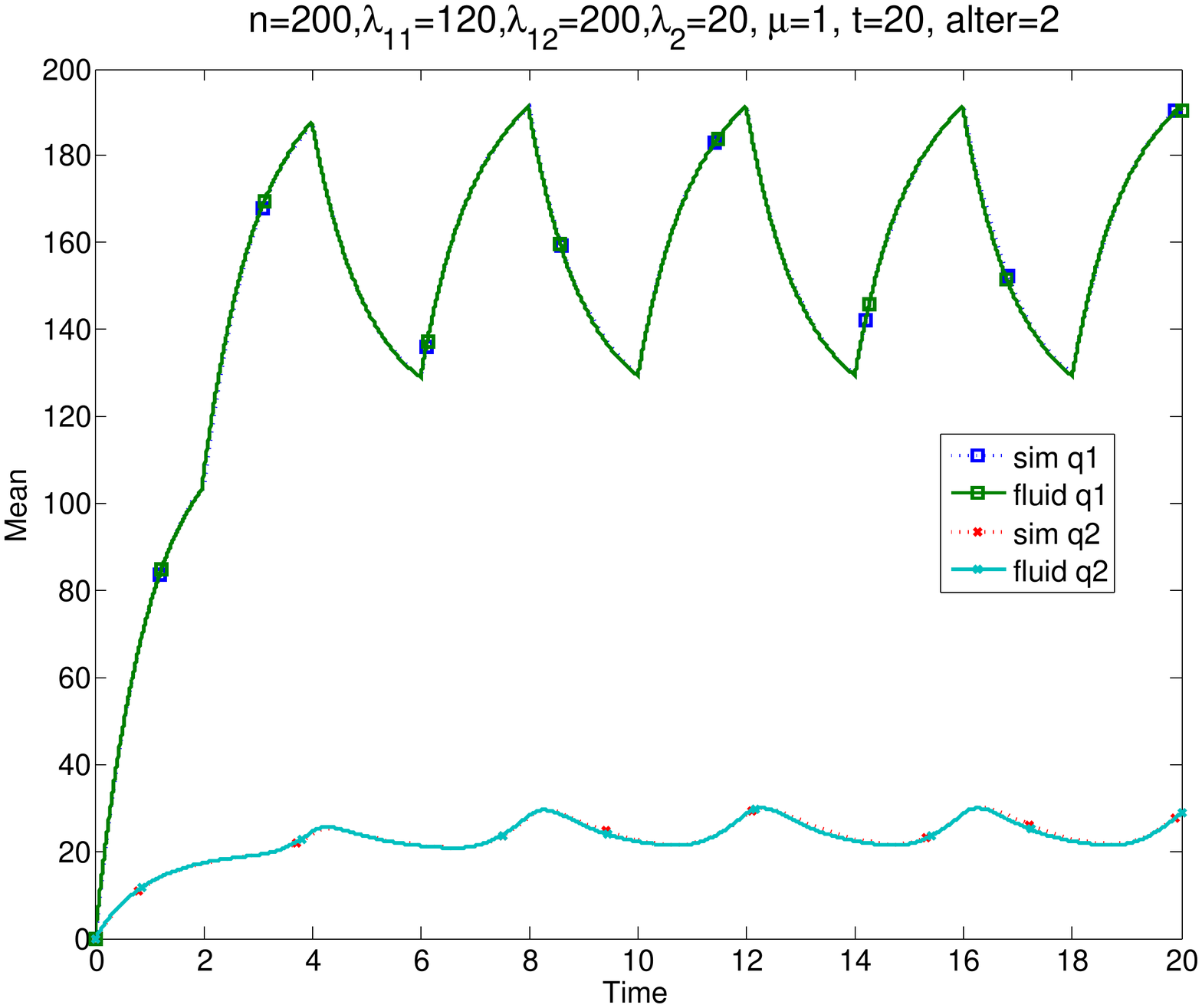}}
    \caption{Comparison of mean values, $E\big[X(t)\big]$}
    \label{fig_multiclass_mean}
  \end{center}
\end{figure}

\begin{figure}[htbp]
  \begin{center}
      \subfigure[ Covariance matrix by assuming measure zero]{
		\includegraphics[width = .48\textwidth]{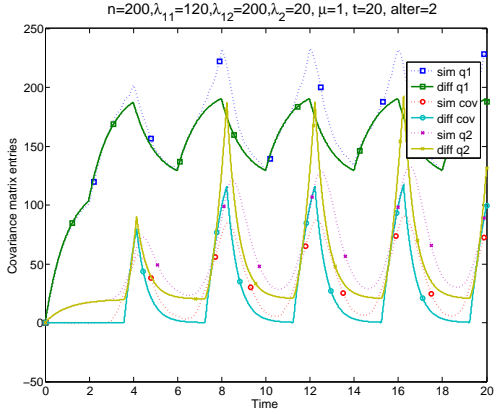}}
      \subfigure[ Covariance matrix by our proposed method]{
		\includegraphics[width = .48\textwidth]{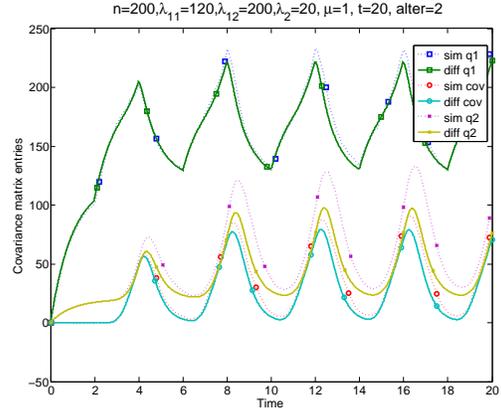}}
    \caption{Comparison of covariance matrix entries, , $Cov\big[X(t),X(t)\big]$}
    \label{fig_multiclass_covariance}
  \end{center}
\end{figure}
Figures \ref{fig_multiclass_mean} and \ref{fig_multiclass_covariance} show the comparison between the method assuming measure zero and our proposed method against simulation. We see that both methods work well for the mean value of $x_1(t)$. However, though not immediately obvious from Figure \ref{fig_multiclass_mean}, there is 5-15\% difference for the mean value of $x_2(t)$ when using the method assuming measure zero while our proposed method shows great accuracy. For the covariance matrix, our proposed method outperforms the method assuming measure zero as seen in Figure \ref{fig_multiclass_covariance}. However, underestimation of variance of $x_2(t)$ is observed in our proposed method. Here, we explain our conjecture on the underestimation of variance. As described in Section \ref{sec_heuristic}, we utilize Gaussian density to obtain new rate functions, $g_i(\cdot,\cdot)$'s. 
In this example, we observe from our numerical experiments that empirical density is not close to Gaussian density when the fluid limit stays near a non-differentiable point. Our conjecture is that asymmetry of empirical density, unlike Gaussian density, causes larger values of covariance matrix entries. However, note that although it does affect the estimation of the covariance matrix (usually underestimation), it does not affect the estimation of mean values of the system state significantly, i.e. we still have the accurate estimation of mean values.
\subsection{Peer networks} \label{subsec_p2p}
\begin{figure}
\centering
\includegraphics[width = .5\textwidth]{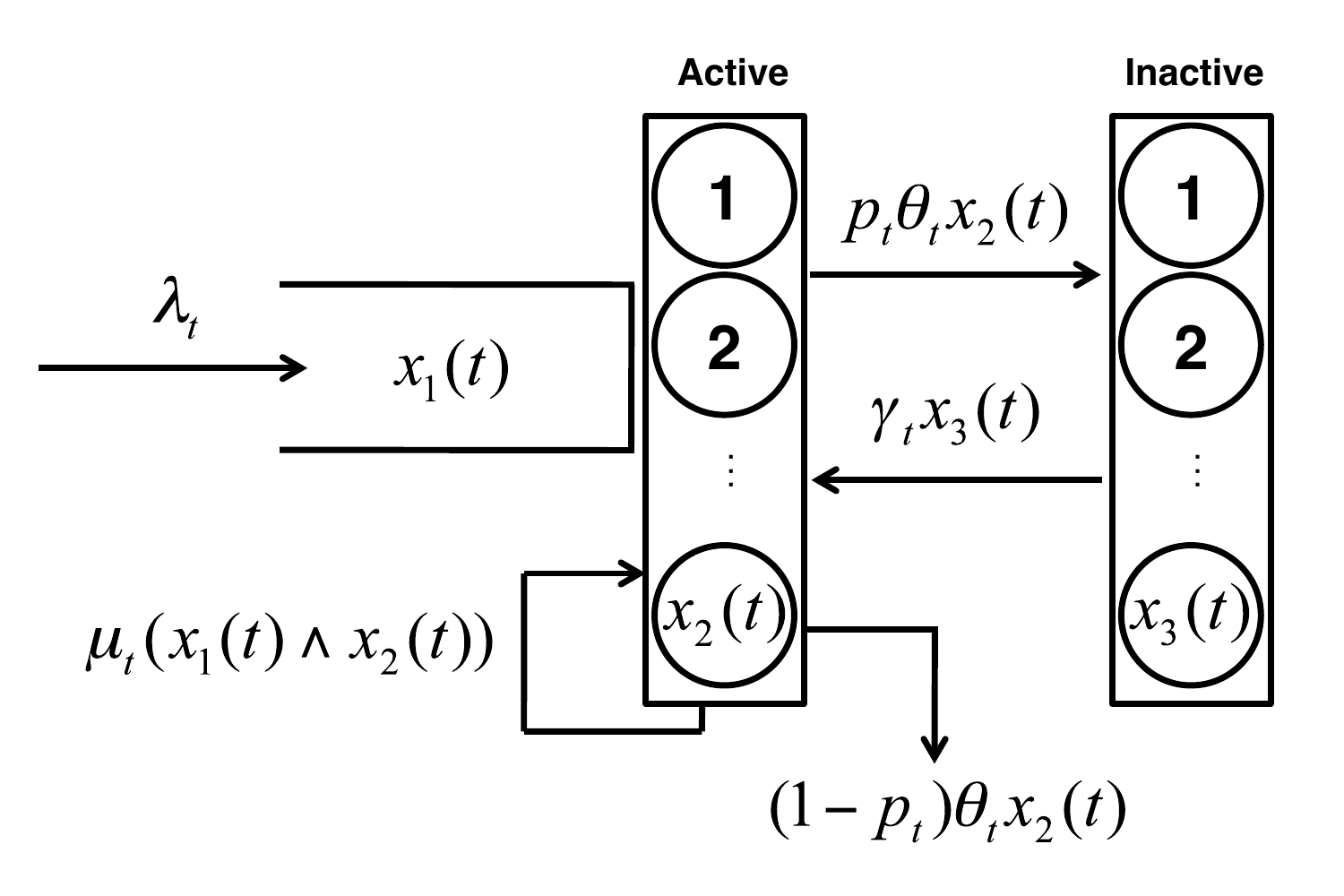}
\caption{Peer networks} \label{fig_p2p}
\end{figure}
Figure \ref{fig_p2p} illustrates the queueing system we consider in this section. Motivated by peer networks with centralized controller (frequently encountered in multimedia content delivering industry), we consider the following scenario of a queueing system.\\
At time $t$, customers arrive to the system at rate $\lambda_t$. Customers in the queue are sent to available active servers. The service rate of each server is $\mu_t$. After a customer is fully served, s/he becomes a new active server. Each server serves customers for a random amount of time with mean $1/\theta_t$, and then either becomes inactive with probability $p_t$ or leaves the system with probability $1-p_t$. Each inactive server spends a random amount of time with mean $1/\gamma_t$ and returns to be active. Note that only active servers can serve customers.\\
Let $X(t)=\big(x_1(t),x_2(t),x_3(t)\big)$ be the state of the system at time $t$ where $x_1(t)$, $x_2(t)$, and $x_3(t)$ are the number of customers, active servers, and inactive servers respectively. Then, $X(t)$ is obtained by solving the following integral equation:
\begin{eqnarray*}
  x_1(t) &=& Y_1\Big(\int_{0}^{t}\lambda_s ds\Big) - Y_2\Big(\int_{0}^{t} \mu_s (x_1(s) \wedge x_2(s)) ds\Big), \\
  x_2(t) &=& Y_2\Big(\int_{0}^{t} \mu_s (x_1(s) \wedge x_2(s)) ds\Big) - Y_3\Big(\int_{0}^{t} \theta_s p_s x_2(s) ds\Big) \\
  &&- Y_4\Big(\int_{0}^{t} \theta_s (1-p_s) x_2(s) ds\Big) + Y_5\Big(\int_{0}^{t} \gamma_s x_3(s) ds\Big), \textrm{ and} \\
  x_3(t) &=& Y_3\Big(\int_{0}^{t} \theta_s p_s x_2(s) ds\Big) - Y_5\Big(\int_{0}^{t} \gamma_s x_3(s) ds\Big),
\end{eqnarray*}
where $Y_1(\cdot)$, $Y_2(\cdot)$, $Y_3(\cdot)$, $Y_4(\cdot)$, and $Y_5(\cdot)$ are independent Poisson processes with rate 1 corresponding to customer arrival, service completion, server's going up, server's leaving, and server's going down respectively. Similar settings are found in the previous research studies (\citet{qiu02}, and \citet{yang04}). They, however, used Poisson processes (i.e. constant rate functions) to construct the system model and did not consider the time-varying properties. Furthermore, they mainly focused on the steady state analysis not providing in-depth analysis of the transient behavior of the system. Therefore, for this application, we can say that we provide more general settings than previous research in that we address time-varying rate functions and provide performance measures through the entire lifespan of the system. In this section, however, we provide the analysis of transient period when the system does not have enough servers to serve customers. This transient analysis is important since the issues on quality of service may arise during this period. We will provide details of the analysis of peer networks in a forthcoming paper.\\
\begin{figure}[htbp]
  \begin{center}
      \subfigure[ Mean numbers by assuming measure zero]{
		\includegraphics[width = .48\textwidth]{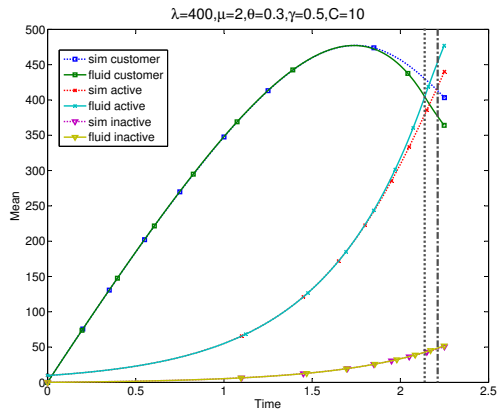}}
      \subfigure[ Mean numbers by our proposed model]{
		\includegraphics[width = .48\textwidth]{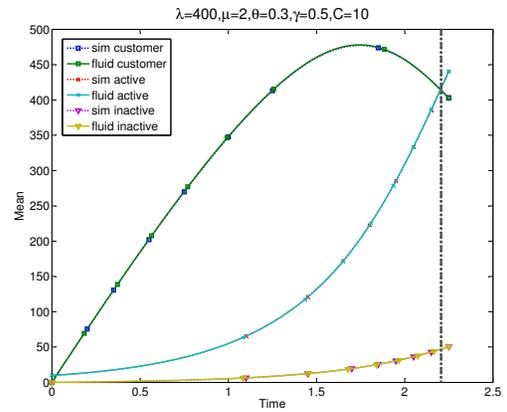}}
    \caption{Comparison of mean values, $E\big[X(t)\big]$}
    \label{fig_p2p_mean}
  \end{center}
\end{figure}
\begin{figure}[htbp]
  \begin{center}
      \subfigure[ Covariance matrix by assuming measure zero]{
		\includegraphics[width = .48\textwidth]{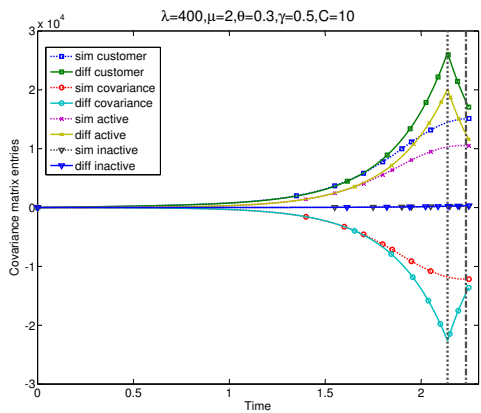}}
      \subfigure[ Covariance matrix by our proposed method]{
		\includegraphics[width = .48\textwidth]{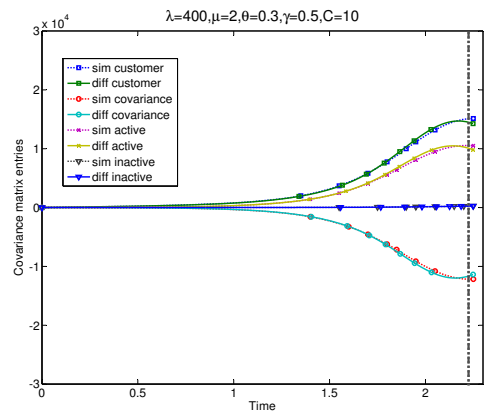}}
    \caption{Comparison of covariance matrix entries, $Cov\big[X(t),X(t)\big]$}
    \label{fig_p2p_covariance}
  \end{center}
\end{figure}
Figure \ref{fig_p2p_mean} illustrates mean values of customers ($E\big[x_1(t)\big]$), active servers ($E\big[x_2(t)\big]$), and inactive servers ($E\big[x_3(t)\big]$) over time. We apply both methods (i.e. the method assuming measure zero and our proposed method) and compare them with the simulation result. We use parameters with $\lambda_t = 400$, $\mu_t= 2$, $\theta_t = 0.3$, $\gamma_t=0.5$, $p_t=0.9$ for $t\ge 0$ as well as $x_1(0)=x_3(0)=0$ and $x_2(0)=10$. We conduct 5,000 simulation runs and obtain mean values by averaging them. In this numerical example, we want to see what happens when the fluid limit becomes close to the critically loaded time point. Figure \ref{fig_p2p_mean} (a) shows comparison between simulation and the method assuming measure zero. We see that it works well when the fluid limits are far from the critically loaded time point. However, as it becomes closer to that point, the method assuming measure zero shows difference from the simulation result. On the other hand, when we apply our proposed method, as seen in Figure \ref{fig_p2p_mean} (b), it provides the almost exact estimation even if the fluid limits are close to the critically loaded point. Figure \ref{fig_p2p_covariance} is the graph of the covariance matrix entries of the system over time. We observe sharp spikes in the method assuming measure zero (Figure \ref{fig_p2p_covariance} (a)). Note that the time when the extreme points of sharp spikes occur is exactly same as the time when the fluid limit hits the critically loaded time point in Figure \ref{fig_p2p_mean} (a) which is the non-differentiable point. However, when we apply our proposed method, we have no sharp spikes at all and the covariance matrix entries are quite close to the simulation result as seen in Figure \ref{fig_p2p_covariance} (b). Thus we believe that our proposed method works well even under this complex scenario.

\section{Conclusion} \label{sec_conclusion}
In this paper, we initially explain the fluid and diffusion models used in analysis of state-dependent queues and show potential problems that one faces in balancing accuracy and computational tractability. The first problem is from the fact that expectation of a function of a random vector $X$ is not equal to the value of the function of the expectation of $X$. Therefore, unless they are equal or close, the fluid model may not provide an accurate estimation of mean values of the system state. The second problem is caused by non-differentiability of rate functions which prevents applying the diffusion model in \citet{kurtz78}. Therefore, addressing these problems is quite important in order to develop accurate approximations as well as to achieve computational feasibility. For that, we proposed a methodology to obtain the exact estimation of mean values of system states and an algorithm to achieve computational tractability.\\
The basic idea of our approach is to construct a new stochastic process which has the fluid limit exactly same as the mean value of the system state. We proved that if rate functions in the original model satisfy the conditions to apply the fluid model, rate functions in the constructed model also satisfy those conditions. Therefore, we can apply the adjusted fluid model if we can apply the existing fluid model. It turns out that there is, in general, no computational method to obtain the adjusted fluid model exactly and hence we utilize Gaussian density to approximate it. By using Gaussian density, we see that rate functions in the constructed model are smooth and we are able to apply the diffusion model in \citet{kurtz78} even if we could not apply it to the original process.\\
To validate our proposed method, we provide several numerical examples of non-stationary state-dependent queueing systems. In the examples, we observe that our proposed method shows great accuracy compared with the method assuming measure zero (which is the only other way in the literature, to the best of our knowledge, that provides computational tractability). Due to space restriction, we have not shown all examples where our method works well. We, however, observe that when the Gaussian density assumption is inaccurate, especially near non-smooth points, our methodology needs further investigation for the covariance matrix. We conjecture that this phenomenon is from the gap between empirical and Gaussian density. To address this, one can investigate the properties of specific rate functions that affect the shape of empirical density and can devise a new algorithm finding $g_i(\cdot,\cdot)$'s from other density functions in the future.  
\appendix
\appendixpage
\addappheadtotoc
\section{Derivation of $g_i(t,x)$'s} \label{app_gi}
For fixed $t_0>0$, let $n=n_{t_0}$, $\mu_1 = \mu_{t_0}^1$, $\beta = \beta_{t_0}$, $p=p_{t_0}$, $x_1 = x_1(t_0) \sim N(z_1, \sigma_1^2)$, and $x_2 = x_2(t_0) \sim N(z_2, \sigma_2^2)$. For $z=(z_1,z_2)'$, we have
\begin{eqnarray}
	g_3\big(t_0, z\big) &=& E\big[\mu_1(x_1 \wedge n)\big] = \mu_1\Big\{E[x_1\mathbb{I}_{x_1 \le n}]+nPr[x_1 > n]\Big\} \nonumber  \\
		&=& \mu_1\Bigg[\int_{-\infty}^{n} \frac{x}{\sqrt{2\pi}\sigma_1}\exp\bigg(-\frac{(x-z_1)^2}{2\sigma_1^2}\bigg) dx + nPr[x_1 > n]\Bigg] \nonumber \\
                &=& \mu_1\Bigg[\frac{-\sigma_1}{\sqrt{2\pi}}\int_{-\infty}^{n} -\frac{x-z_1}{\sigma_1^2}\exp\bigg(-\frac{(x-z_1)^2}{2\sigma_1^2}\bigg) dx + z_1Pr[x_1 \le n] + nPr[x_1 > n]\Bigg] \nonumber \\
                &=& \mu_1 \Bigg[-\sigma_1^2 \frac{1}{\sqrt{2\pi}\sigma_1} \exp\bigg(-\frac{(n-z_1)^2}{2\sigma_1^2}\bigg)+(z_1-n)Pr(x_1\le n) + n\Bigg]. \nonumber
\end{eqnarray}
Therefore, by making $t_0>0$ arbitrary, we have $g_3(t,x)$.\\
Note $g_4(\cdot,\cdot)$ and $g_5(\cdot,\cdot)$ are same except a constant part with respect to $x$. Therefore, it is enough to derive $g_5(\cdot,\cdot)$. We can show that
\begin{eqnarray}
	g_5\big(t_0, z\big) &=& E\big[\beta p (x_1 - n)^+\big] = \beta p \big\{E[x_1\vee n]-n\big\} \nonumber \\
        &=& \beta p \Big\{E[x_1\mathbb{I}_{x_1 > n}]+nPr[x_1 \le n] -n\Big\} \nonumber  \\
		&=& \beta p \Bigg[\int_{n}^{\infty} \frac{x}{\sqrt{2\pi}\sigma_1}\exp\bigg(-\frac{(x-z_1)^2}{2\sigma_1^2}\bigg) dx + nPr[x_1 \le n] -n\Bigg] \nonumber \\
                &=& \beta p \Bigg[\frac{-\sigma_1}{\sqrt{2\pi}}\int_{n}^{\infty} -\frac{x-z_1}{\sigma_1^2}\exp\bigg(-\frac{(x-z_1)^2}{2\sigma_1^2}\bigg) dx + z_1Pr[x_1 > n] + nPr[x_1 \le n] -n \Bigg] \nonumber \\
                &=& \beta p \Bigg[\sigma_1^2 \frac{1}{\sqrt{2\pi}\sigma_1} \exp\bigg(-\frac{(n-z_1)^2}{2\sigma_1^2}\bigg)+(z_1-n)Pr(x_1 > n) \Bigg]. \nonumber
\end{eqnarray}
Therefore, by making $t_0>0$ arbitrary, we have $g_5(t,x)$.
\section{Numerical results for Section \ref{sec_application}} \label{app_table}
\begin{table}[htdp]
\caption{Estimation of $E\big[x_1(t)\big]$ over time; difference from simulation}
\begin{center}
\footnotesize
\begin{tabular}{|c|c|c|c|c|c|c|c|c|c|c|c|}
\hline
\multicolumn{2}{|c|}{Experiments} & \multicolumn{10}{|c|}{Time ($t$)} \\ \hline 
\# & type & 6 & 7 & 8 & 9 & 10 & 11 & 12 & 13 & 14 & 15 \\ \hline
\multirow{2}{*}{1} & proposed & 6.52 & 0.98 & -3.39 & -1.07 & -3.05 & -0.40 & 0.91 & 0.25 & -0.69 & -0.01 \\ \cline{2-12}
& meas. 0 &4.42 & 0.82 & -3.63 & -1.94 & -3.60 & -0.23 & 0.75 & -0.15 & -2.59 & 0.11 \\ \hline
\multirow{2}{*}{2} & proposed & 2.69 & 0.44 & -3.13 & -0.82 & -1.08 & -0.32 & 0.48 & 0.15 & -0.46 & -0.05 \\ \cline{2-12}
& meas. 0 &3.35 & -0.42 & -2.92 & -1.64 & -1.18 & -1.01 & 0.85 & -0.44 & -0.36 & -0.60 \\ \hline
\multirow{2}{*}{3} & proposed & 2.33 & 0.28 & -3.11 & -1.01 & -1.36 & -0.39 & 0.10 & -0.02 & -1.68 & -0.15 \\ \cline{2-12}
& meas. 0 &2.34 & -0.42 & -2.67 & -1.55 & -1.49 & -1.00 & 0.52 & -0.49 & -0.54 & -0.53 \\ \hline
\multirow{2}{*}{4} & proposed & 1.18 & 0.14 & -1.54 & -0.30 & -0.01 & 0.12 & 0.22 & 0.22 & -0.10 & -0.02 \\ \cline{2-12}
& meas. 0 &0.65 & -0.96 & -1.98 & -1.32 & -0.94 & -0.95 & 0.04 & -0.64 & -0.61 & -0.94 \\ \hline
\multirow{2}{*}{5} & proposed & 7.04 & 1.36 & -3.67 & -0.69 & -1.38 & -0.57 & 0.80 & 0.23 & -2.82 & -0.63 \\ \cline{2-12}
& meas. 0 &5.55 & 1.04 & -3.20 & -0.93 & -1.31 & -0.53 & 0.46 & 0.06 & -1.22 & -0.18 \\ \hline
\multirow{2}{*}{6} & proposed & 3.61 & 0.76 & -3.05 & -1.13 & -0.67 & 0.18 & 1.12 & 0.20 & -0.95 & -0.25 \\ \cline{2-12}
& meas. 0 &2.53 & -0.07 & -3.01 & -1.72 & -1.46 & -0.43 & 0.60 & -0.47 & -1.57 & -0.80 \\ \hline
\multirow{2}{*}{7} & proposed & 1.93 & 0.65 & -1.06 & -0.25 & -0.63 & 0.17 & 0.12 & -0.21 & -0.65 & -0.20 \\ \cline{2-12}
& meas. 0 &0.50 & -0.86 & -2.07 & -1.51 & -1.04 & -0.73 & -0.47 & -1.07 & -0.63 & -0.76 \\ \hline
\multirow{2}{*}{8} & proposed & 0.72 & 0.07 & -0.46 & 0.04 & -0.04 & -0.14 & 0.42 & -0.07 & -0.48 & -0.01 \\ \cline{2-12}
& meas. 0 &0.04 & -0.98 & -1.40 & -0.91 & -0.57 & -0.85 & -0.13 & -0.69 & -0.73 & -0.46 \\ \hline
\multirow{2}{*}{9} & proposed & 0.81 & 0.25 & -0.96 & -0.25 & -0.11 & -0.09 & 0.38 & -0.06 & -0.24 & -0.02 \\ \cline{2-12}
& meas. 0 &0.53 & -0.50 & -1.31 & -0.88 & -0.34 & -0.61 & 0.17 & -0.51 & -0.06 & -0.32 \\ \hline
\multirow{2}{*}{10} & proposed & 6.44 & 1.18 & -4.73 & -1.73 & -2.21 & -0.45 & 0.30 & -0.01 & -1.10 & -0.11 \\ \cline{2-12}
& meas. 0 &6.46 & 0.77 & -3.83 & -1.62 & -2.84 & -0.83 & 0.84 & 0.00 & -2.77 & -0.60 \\
\hline
\end{tabular}
\end{center}
\normalsize
\label{tab_mean_x1}
\end{table}

\begin{table}[htdp]
\caption{Estimation of $E\big[x_2(t)\big]$ over time; difference from simulation}
\begin{center}
\footnotesize
\begin{tabular}{|c|c|c|c|c|c|c|c|c|c|c|c|}
\hline
\multicolumn{2}{|c|}{Experiments} & \multicolumn{10}{|c|}{Time ($t$)} \\ \hline 
\# & type & 6 & 7 & 8 & 9 & 10 & 11 & 12 & 13 & 14 & 15 \\ \hline
\multirow{2}{*}{1} & proposed & -2.00 & 3.50 & 2.36 & -0.53 & 0.57 & -1.00 & -0.99 & -0.30 & -0.44 & -0.76 \\ \cline{2-12}
& meas. 0 &11.68 & 12.60 & 7.38 & 5.88 & 11.64 & 8.18 & 5.24 & 6.29 & 10.47 & 7.82 \\ \hline
\multirow{2}{*}{2} & proposed & -2.22 & 2.71 & 1.90 & -2.44 & -0.94 & -1.82 & -0.91 & -0.10 & -0.38 & -0.76 \\ \cline{2-12}
& meas. 0 &45.00 & 53.07 & 33.49 & 37.12 & 41.73 & 44.51 & 31.55 & 37.48 & 40.57 & 43.21 \\ \hline
\multirow{2}{*}{3} & proposed & -2.49 & 1.88 & 1.00 & -3.58 & -2.09 & -3.32 & -3.08 & -3.01 & -3.15 & -4.02 \\ \cline{2-12}
& meas. 0 &28.64 & 37.65 & 19.44 & 21.88 & 24.73 & 28.38 & 16.37 & 21.45 & 22.77 & 26.96 \\ \hline
\multirow{2}{*}{4} & proposed & 0.24 & 2.66 & 1.35 & -1.68 & -0.91 & -0.19 & 0.25 & 0.45 & 0.02 & -0.53 \\ \cline{2-12}
& meas. 0 &67.95 & 69.81 & 47.03 & 51.81 & 56.69 & 59.66 & 45.16 & 50.75 & 54.48 & 57.27 \\ \hline
\multirow{2}{*}{5} & proposed & -1.01 & 4.41 & 3.16 & -0.05 & 1.55 & 0.57 & -0.58 & -0.07 & -0.09 & -1.63 \\ \cline{2-12}
& meas. 0 &9.61 & 12.28 & 7.50 & 5.73 & 11.28 & 9.42 & 5.36 & 6.00 & 9.51 & 8.02 \\ \hline
\multirow{2}{*}{6} & proposed & -2.63 & 2.48 & 2.23 & -2.39 & -1.32 & -0.84 & -0.12 & 1.04 & 0.89 & -0.04 \\ \cline{2-12}
& meas. 0 &44.23 & 51.84 & 32.45 & 35.11 & 39.27 & 43.72 & 31.00 & 35.94 & 38.83 & 41.83 \\ \hline
\multirow{2}{*}{7} & proposed & 0.33 & 3.42 & 3.00 & 1.01 & 0.70 & 0.25 & 0.41 & 0.71 & 0.27 & -0.17 \\ \cline{2-12}
& meas. 0 &78.08 & 78.96 & 60.84 & 64.86 & 69.59 & 71.19 & 59.15 & 63.20 & 67.42 & 68.95 \\ \hline
\multirow{2}{*}{8} & proposed & 2.81 & 3.03 & 2.40 & 1.45 & 1.29 & 0.58 & -0.08 & 0.41 & 0.12 & -1.11 \\ \cline{2-12}
& meas. 0 &92.68 & 90.60 & 73.97 & 77.06 & 80.98 & 81.48 & 70.82 & 74.24 & 77.96 & 78.55 \\ \hline
\multirow{2}{*}{9} & proposed & -0.86 & 1.25 & 1.44 & -0.77 & -0.19 & -0.18 & 0.08 & 0.84 & 0.42 & 0.41 \\ \cline{2-12}
& meas. 0 &80.15 & 79.90 & 57.59 & 62.03 & 67.09 & 68.91 & 55.09 & 59.98 & 64.19 & 66.14 \\ \hline
\multirow{2}{*}{10} & proposed & -2.67 & 6.62 & 3.79 & -2.50 & -0.49 & -2.18 & -1.99 & -1.35 & -1.38 & -1.21 \\ \cline{2-12}
& meas. 0 &8.53 & 23.91 & 10.73 & 8.77 & 10.78 & 13.05 & 5.67 & 8.96 & 9.13 & 11.69 \\
\hline
\end{tabular}
\end{center}
\normalsize
\label{tab_mean_x2}
\end{table}

\begin{table}[htdp]
\caption{Estimation of $Var\big[x_1(t)\big]$ over time; difference from simulation}
\begin{center}
\footnotesize
\begin{tabular}{|c|c|c|c|c|c|c|c|c|c|c|c|}
\hline
\multicolumn{2}{|c|}{Experiments} & \multicolumn{10}{|c|}{Time ($t$)} \\ \hline 
\# & method & 6 & 7 & 8 & 9 & 10 & 11 & 12 & 13 & 14 & 15 \\ \hline
\multirow{2}{*}{1} & proposed & 6.94 & 0.94 & -1.92 & -2.02 & -3.66 & -0.20 & 1.70 & -1.02 & 0.49 & 2.89 \\ \cline{2-12}
& meas. 0 &-11.03 & 2.93 & -1.93 & 17.31 & -24.44 & 1.42 & 1.66 & 14.89 & -19.73 & 4.16 \\ \hline
\multirow{2}{*}{2} & proposed & 2.84 & 3.83 & -6.05 & -0.10 & -0.50 & 4.24 & 1.62 & 2.67 & -1.02 & 1.62 \\ \cline{2-12}
& meas. 0 &-6.28 & 16.69 & 6.76 & -14.45 & -12.97 & 17.62 & 12.90 & -11.61 & -14.51 & 15.29 \\ \hline
\multirow{2}{*}{3} & proposed & 4.15 & 2.09 & -0.60 & 2.15 & -6.57 & 0.50 & -3.10 & 1.15 & 2.76 & 3.74 \\ \cline{2-12}
& meas. 0 &-0.56 & 13.38 & 7.30 & -8.74 & -12.97 & 11.92 & 4.53 & -10.16 & -2.13 & 14.22 \\ \hline
\multirow{2}{*}{4} & proposed & -0.52 & -4.36 & -2.81 & 3.07 & -0.03 & 2.96 & 0.79 & 1.27 & 3.30 & 0.35 \\ \cline{2-12}
& meas. 0 &-16.38 & 11.18 & 14.13 & -12.86 & -17.81 & 17.93 & 16.94 & -15.33 & -14.05 & 15.32 \\ \hline
\multirow{2}{*}{5} & proposed & 6.83 & -0.22 & -2.49 & 0.09 & -1.67 & -3.27 & 1.71 & -4.14 & -0.55 & 1.98 \\ \cline{2-12}
& meas. 0 &-2.30 & 1.03 & -1.69 & 10.07 & -10.59 & -1.97 & 1.43 & 5.09 & -7.69 & 3.42 \\ \hline
\multirow{2}{*}{6} & proposed & 5.22 & 0.62 & -6.25 & -0.81 & -4.32 & -1.95 & 4.41 & 1.97 & -0.61 & 4.93 \\ \cline{2-12}
& meas. 0 &-1.19 & 7.72 & 1.39 & -7.42 & -11.70 & 5.61 & 10.28 & -4.33 & -7.73 & 12.15 \\ \hline
\multirow{2}{*}{7} & proposed & 2.91 & -2.29 & -1.04 & 0.92 & 0.21 & 0.18 & 3.14 & -1.10 & 4.36 & 2.28 \\ \cline{2-12}
& meas. 0 &-17.83 & 14.52 & 18.27 & -16.55 & -22.37 & 17.07 & 20.88 & -18.77 & -18.14 & 19.01 \\ \hline
\multirow{2}{*}{8} & proposed & -1.79 & 0.65 & -0.43 & 0.83 & 3.35 & -0.71 & 3.63 & 2.10 & 1.85 & 0.72 \\ \cline{2-12}
& meas. 0 &-26.38 & 16.44 & 21.26 & -18.37 & -22.66 & 16.73 & 23.72 & -17.42 & -25.80 & 18.25 \\ \hline
\multirow{2}{*}{9} & proposed & 0.62 & -0.86 & -0.83 & 3.53 & 3.36 & 5.09 & 1.52 & 1.71 & 2.73 & -1.37 \\ \cline{2-12}
& meas. 0 &-17.84 & 13.72 & 17.09 & -14.12 & -17.40 & 19.78 & 18.07 & -16.57 & -19.03 & 14.36 \\ \hline
\multirow{2}{*}{10} & proposed & 4.48 & -0.32 & -9.84 & 1.26 & -4.24 & 3.37 & 1.32 & 1.00 & 0.22 & 1.15 \\ \cline{2-12}
& meas. 0 &4.12 & 7.27 & -6.22 & -2.69 & -5.55 & 10.87 & 3.68 & -3.27 & -2.12 & 8.98 \\
\hline
\end{tabular}
\end{center}
\normalsize
\label{tab_var_x1}
\end{table}

\begin{table}[htdp]
\caption{Estimation of $Cov\big[x_1(t),x_2(t)\big]$ over time; difference from simulation}
\begin{center}
\footnotesize
\begin{tabular}{|c|c|c|c|c|c|c|c|c|c|c|c|}
\hline
\multicolumn{2}{|c|}{Experiments} & \multicolumn{10}{|c|}{Time ($t$)} \\ \hline 
\# & type & 6 & 7 & 8 & 9 & 10 & 11 & 12 & 13 & 14 & 15 \\ \hline
\multirow{2}{*}{1} & proposed & -3.03 & -3.27 & 4.75 & 3.10 & -1.72 & -3.63 & 0.39 & -4.00 & -4.15 & 0.91 \\ \cline{2-12}
& meas. 0 &25.05 & -3.88 & 4.43 & -6.75 & 15.78 & -4.69 & -0.30 & -11.23 & 4.26 & -2.03 \\ \hline
\multirow{2}{*}{2} & proposed & -6.76 & 7.23 & 2.87 & -4.73 & 11.50 & -0.47 & 6.09 & 5.48 & 5.86 & 4.82 \\ \cline{2-12}
& meas. 0 &29.60 & -6.12 & -6.56 & 3.81 & 36.90 & -13.05 & -1.41 & 12.63 & 31.69 & -5.53 \\ \hline
\multirow{2}{*}{3} & proposed & -6.74 & -2.24 & 4.53 & -9.51 & -28.97 & -3.44 & -2.57 & -6.43 & 5.52 & -3.76 \\ \cline{2-12}
& meas. 0 &25.67 & -15.55 & 0.47 & 6.26 & 6.46 & -15.26 & -6.35 & 8.68 & 31.03 & -13.44 \\ \hline
\multirow{2}{*}{4} & proposed & -0.01 & -13.57 & -8.53 & -12.42 & -9.73 & -0.00 & -10.28 & -16.74 & -1.43 & -11.64 \\ \cline{2-12}
& meas. 0 &58.61 & -29.39 & -21.29 & 10.14 & 44.87 & -14.34 & -22.28 & 5.51 & 46.98 & -26.05 \\ \hline
\multirow{2}{*}{5} & proposed & -7.19 & -0.18 & 0.13 & 2.88 & 7.13 & -8.20 & -0.42 & 1.07 & 4.17 & -1.88 \\ \cline{2-12}
& meas. 0 &19.73 & -4.03 & -1.03 & -12.97 & 26.24 & -11.11 & -1.32 & -12.31 & 20.40 & -4.25 \\ \hline
\multirow{2}{*}{6} & proposed & 2.91 & 4.90 & 2.63 & -7.64 & -7.99 & 1.17 & 2.14 & 3.80 & 6.94 & 7.92 \\ \cline{2-12}
& meas. 0 &37.93 & -15.54 & -15.96 & -3.32 & 25.88 & -18.86 & -15.00 & 6.43 & 34.79 & -10.85 \\ \hline
\multirow{2}{*}{7} & proposed & -4.38 & -1.88 & 0.97 & -14.36 & 3.15 & -1.95 & -1.19 & -1.08 & -4.77 & -0.46 \\ \cline{2-12}
& meas. 0 &52.91 & -16.88 & -16.53 & -3.19 & 43.26 & -15.54 & -16.02 & 6.73 & 34.99 & -12.29 \\ \hline
\multirow{2}{*}{8} & proposed & -20.99 & -4.21 & -6.33 & -4.51 & 3.12 & -3.43 & -1.85 & -6.78 & -1.81 & -0.79 \\ \cline{2-12}
& meas. 0 &64.94 & -8.85 & -22.74 & 13.30 & 51.79 & -11.89 & -15.66 & 7.80 & 44.16 & -8.72 \\ \hline
\multirow{2}{*}{9} & proposed & -15.01 & -6.15 & -6.27 & 3.33 & -0.25 & 2.45 & -6.00 & -7.57 & -6.93 & -6.74 \\ \cline{2-12}
& meas. 0 &55.84 & -12.34 & -17.97 & 19.55 & 45.76 & -4.17 & -15.42 & 7.67 & 37.97 & -12.70 \\ \hline
\multirow{2}{*}{10} & proposed & -18.70 & 7.57 & -3.70 & -4.76 & 8.09 & -6.43 & -2.86 & -0.03 & 2.95 & -1.11 \\ \cline{2-12}
& meas. 0 &-21.43 & -2.63 & -5.67 & -0.67 & 8.99 & -15.71 & -4.40 & 3.66 & 4.18 & -9.87 \\
\hline
\end{tabular}
\end{center}
\normalsize
\label{tab_cov}
\end{table}

\begin{table}[htdp]
\caption{Estimation of $Var\big[x_2(t)\big]$ over time; difference from simulation}
\begin{center}
\footnotesize
\begin{tabular}{|c|c|c|c|c|c|c|c|c|c|c|c|}
\hline
\multicolumn{2}{|c|}{Experiments} & \multicolumn{10}{|c|}{Time ($t$)} \\ \hline 
\# & type & 6 & 7 & 8 & 9 & 10 & 11 & 12 & 13 & 14 & 15 \\ \hline
\multirow{2}{*}{1} & proposed & -2.15 & 3.52 & 1.48 & -0.72 & -0.34 & -0.45 & 0.78 & 1.59 & 1.31 & 0.83 \\ \cline{2-12}
& meas. 0 &6.74 & 5.31 & 2.34 & -8.01 & 5.46 & 1.00 & 1.84 & -2.78 & 2.81 & -1.20 \\ \hline
\multirow{2}{*}{2} & proposed & 1.29 & 9.81 & 8.50 & 3.31 & 7.72 & 6.70 & 6.05 & 5.84 & 5.42 & 5.91 \\ \cline{2-12}
& meas. 0 &7.06 & 14.88 & -6.56 & -0.09 & 17.07 & 12.12 & -3.90 & 5.44 & 16.94 & 13.19 \\ \hline
\multirow{2}{*}{3} & proposed & -5.60 & 2.13 & -0.71 & -5.14 & -1.93 & -3.21 & -2.71 & -1.18 & -0.61 & -0.79 \\ \cline{2-12}
& meas. 0 &-2.22 & 4.53 & -10.37 & -5.14 & 4.01 & 1.00 & -8.89 & 0.66 & 6.45 & 5.96 \\ \hline
\multirow{2}{*}{4} & proposed & 5.71 & 8.34 & 2.63 & -2.01 & -0.83 & 2.01 & -0.03 & -0.27 & 0.51 & 2.66 \\ \cline{2-12}
& meas. 0 &28.49 & 26.24 & -13.87 & -3.06 & 15.93 & 15.77 & -12.78 & 0.50 & 17.18 & 17.62 \\ \hline
\multirow{2}{*}{5} & proposed & -0.97 & 4.30 & 1.25 & -0.07 & 3.22 & 3.50 & -0.22 & -0.16 & 1.50 & 0.76 \\ \cline{2-12}
& meas. 0 &3.67 & 5.28 & 1.60 & -4.35 & 7.92 & 6.18 & 1.65 & -2.36 & 5.70 & 4.03 \\ \hline
\multirow{2}{*}{6} & proposed & 2.23 & 11.10 & 8.33 & 2.94 & 4.21 & 4.30 & 1.38 & 2.79 & 3.63 & 5.03 \\ \cline{2-12}
& meas. 0 &11.13 & 21.61 & -3.19 & -0.25 & 12.83 & 14.35 & -6.15 & 1.31 & 12.95 & 14.75 \\ \hline
\multirow{2}{*}{7} & proposed & 5.23 & 7.48 & 5.57 & 1.60 & 2.20 & 4.24 & 3.45 & 4.88 & 5.03 & 4.33 \\ \cline{2-12}
& meas. 0 &33.03 & 27.56 & -16.08 & -4.39 & 21.67 & 19.11 & -11.03 & 2.36 & 24.85 & 19.64 \\ \hline
\multirow{2}{*}{8} & proposed & 10.11 & 7.03 & 3.99 & 2.44 & 3.47 & 2.45 & 2.53 & 2.25 & 1.73 & 2.30 \\ \cline{2-12}
& meas. 0 &62.03 & 46.52 & -13.63 & 0.58 & 30.10 & 23.25 & -13.94 & -0.09 & 26.38 & 20.60 \\ \hline
\multirow{2}{*}{9} & proposed & 8.18 & 7.49 & 3.22 & 0.53 & 3.83 & 4.55 & 3.14 & 4.24 & 4.90 & 5.28 \\ \cline{2-12}
& meas. 0 &39.93 & 31.88 & -18.20 & -4.36 & 22.39 & 18.66 & -12.73 & 1.21 & 23.00 & 18.98 \\ \hline
\multirow{2}{*}{10} & proposed & -0.34 & 12.31 & 5.05 & -2.01 & 1.38 & 1.13 & -3.01 & -3.64 & -4.35 & -1.66 \\ \cline{2-12}
& meas. 0 &-5.73 & 7.15 & -1.91 & -3.68 & 0.93 & -2.52 & -8.00 & -4.34 & -3.94 & -5.72 \\
\hline
\end{tabular}
\end{center}
\normalsize
\label{tab_var_x2}
\end{table}
\newpage
\bibliographystyle{plainnat}
\bibliography{paper}

\begin{thebibliography}{10}

\bibitem{arnold92}
Ludwig Arnold.
\newblock {\em Stochastic Differential Equations: Theory and Applications}.
\newblock Krieger Publishing Company, 1992.

\bibitem{ethier86}
Stewart~N. Ethier and Thomas~G. Kurtz.
\newblock {\em Markov Processes: Characterization and Convergence}.
\newblock A John Wiley \& Sons, Inc., Publication, 1 edition, 1986.

\bibitem{folland99}
Gerald~B. Folland.
\newblock {\em Real Analysis : Modern Techniques and Their Applications}.
\newblock A John Wiley \& Sons, Inc., Publication, 2 edition, 1999.

\bibitem{iglehart65}
Donald~L. Iglehart.
\newblock Limiting diffusion approximations for the many server queue and the
  repairman problem.
\newblock {\em Journal of Applied Probability}, 2(2):429--441, December 1965.

\bibitem{kurtz78}
Thomas~G. Kurtz.
\newblock Strong approximation theorems for density dependent markov chains.
\newblock {\em Stochastic Processes and their Applications}, 6(3):223--240, feb
  1978.

\bibitem{mandelbaum98}
Avi Mandelbaum, William~A. Massey, and Martin~I. Reiman.
\newblock Strong approximations for markovian service networks.
\newblock {\em Queueing Systems}, 30:149--201, 1998.

\bibitem{mandelbaum02}
Avi Mandelbaum, William~A. Massey, and Brian Rider.
\newblock Queue lengths and waiting times for multiserver queues with
  abandonment and retrials.
\newblock {\em Telecommunication Systems}, 21(2-4):149--171, 2002.

\bibitem{mandelbaum98-2}
Avi Mandelbaum and Gennady Pats.
\newblock State-dependent stochastic networks. part i: Approximations and
  applications with continuous diffusion limits.
\newblock {\em The Annals of Applied Probability}, 8(2):569--646, may 1998.

\bibitem{Whitt04}
Ward Whitt.
\newblock Efficiency-driven heavy-traffic approximations for many-server queues
  with abandonments.
\newblock {\em Management Science}, 50(10):1449--1461, oct 2006.

\bibitem{Whitt06b}
Ward Whitt.
\newblock Fluid models for multiserver queues with abandonments.
\newblock {\em Operations Research}, 54(1):37--54, January-February 2006.

\end{thebibliography}

\end{document}